\documentclass[a4paper,UKenglish]{lipics_md}
\usepackage{microtype}
\usepackage{mathtools}
\usepackage[normalem]{ulem}

\usepackage{amsmath}
\usepackage{amsfonts}
\usepackage{arydshln,leftidx,mathtools}
\title{Weighted Poisson--Delaunay Mosaics\footnote{This project has received funding from the
       European Research Council (ERC) under the European Union's Horizon 2020
       research and innovation programme (grant agreement No 78818 Alpha).
       It is also partially supported by the DFG Collaborative Research Center
       TRR 109, `Discretization in Geometry and Dynamics',
       through grant no.\ I02979-N35 of the Austrian Science Fund (FWF).}}
\titlerunning{Weighted Poisson--Delaunay Mosaics}

\author[]{Herbert Edelsbrunner}
\author[]{Anton Nikitenko}
\affil[]{IST Austria (Institute of Science and Technology Austria),
  Am Campus 1, \\ 3400 Klosterneuburg, Austria,
  \texttt{edels@ist.ac.at}, \texttt{anton.nikitenko@ist.ac.at}}
\authorrunning{H. Edelsbrunner and A. Nikitenko}

\mathsubjclass{60D05 Geometric probability and stochastic geometry,
  68U05 Computer graphics; computational geometry.} 

\newcommand{\mm}[1] {\ifmmode{#1}\else{\mbox{\(#1\)}}\fi}

\newcommand{\denselist}{\itemsep 0pt\parsep=1pt\partopsep 0pt}

\newcommand{\ourproof}{\begin{proof}}
\newcommand{\eop}{\end{proof}}  %
\newcommand{\UpperHalfPlane}{\mm{{\mathbf H}}}
\newcommand{\Rspace}        {\mm{{\mathbb R}}}
\newcommand{\Sspace}        {\mm{{\mathbb S}}}
\newcommand{\Sphere}        {\mm{{S}}}
\newcommand{\Expected}[1]   {\mm{{\mathbb E}{[{#1}]}}}
\newcommand{\PP}            {\mm{{\mathbb P}}}
\newcommand{\Probable}[1]   {\mm{{\PP}{[{#1}]}}}
\newcommand{\PemptyOnly}    {\mm{{\PP}_{\emptyset}}}
\newcommand{\Pempty}[1]     {\mm{{\PP}_{\emptyset}{({#1})}}}
\newcommand{\Random}[4]     {\mm{{U}_{{#1},{#2}}^{{#3},{#4}}}}
\newcommand{\One}           {\mm{{{\bf 1}}}}

\newcommand{\Gama}[1]       {\mm{{\Gamma}{\left({#1}\right)}}}
\newcommand{\iGama}[2]      {\mm{{\gamma}{\left({#1};\,{#2}\right)}}}
\newcommand{\Beta}[2]       {\mm{{B}{\left({#1},{#2}\right)}}}
\newcommand{\iBeta}[3]      {\mm{{B}_{#1}{({#2},{#3})}}}
\newcommand{\Hyper}[2]      {\mm{_{#1}{F}_{#2}}}
\newcommand{\HyperReg}[2]   {\mm{_{#1}\tilde{F}_{#2}}}
\newcommand{\LGrass}[2]     {\mm{{\cal L}_{#1}^{#2}}}
\newcommand{\Normal}[2]     {\mm{{\cal N}{({#1},{#2})}}}
\newcommand{\Delaunay}[1]   {\mm{{\rm Del}\,{#1}}}
\newcommand{\Del}[2]   {\mm{{\rm Del}_{#1}({#2})}}
\newcommand{\Voronoi}[1]    {\mm{{\rm Vor}\,{#1}}}
\newcommand{\voronoi}[1]    {\mm{{\rm domain}{({#1})}}}

\newcommand{\Rfun}          {\mm{{\cal R}}}
\newcommand{\density}       {\mm{\rho}}
\newcommand{\intensity}[1]  {\mm{\varrho}{({#1})}}

\newcommand{\ttt}           {\mm{{\bf t}}}
\newcommand{\uuu}           {\mm{{\bf u}}}
\newcommand{\xxx}           {\mm{{\bf x}}}

\newcommand{\Vis}[1]        {\mm{{\rm Vis}{({#1})}}}

\newcommand{\Ccon}[4]       {\mm{{C}_{{#1},{#2}}^{{#3},{#4}}}}
\newcommand{\Dcon}[3]       {\mm{{D}_{#1}^{{#2},{#3}}}}
\newcommand{\ccon}[4]       {\mm{{c}_{{#1},{#2}}^{{#3},{#4}}}}
\newcommand{\dcon}[3]       {\mm{{d}_{#1}^{{#2},{#3}}}}
\newcommand{\Vol}[2]        {\mm{\rm Vol}_{#1}{({#2})}}

\newcommand{\abs}[1]        {\mm{\rm abs\,}{#1}}

\newcommand{\dime}[1]       {\mm{\rm dim\,}{#1}}

\newcommand{\aff}[1]        {\mm{\rm aff\,}{#1}}

\newcommand{\diff}          {\mm{\rm \,d}}
\newcommand{\norm}[1]       {\mm{\|{#1}\|}}
\newcommand{\Edist}[2]      {\mm{\|{#1}-{#2}\|}}

\newcommand{\Region}        {\mm{\Omega}}

\newcommand{\ourparagraph}[1]  {\vspace{0.1in} \noindent \textbf{#1}}
\newcommand{\Skip}[1]       {}

\begin{document}
\maketitle

\begin{abstract}
  Slicing a Voronoi tessellation in $\Rspace^n$ with a $k$-plane gives a
  $k$-dimensional weighted Voronoi tessellation, also known as power diagram
  or Laguerre tessellation.
  Mapping every simplex of the dual weighted Delaunay mosaic to the
  radius of the smallest empty circumscribed sphere whose center lies
  in the $k$-plane gives a generalized discrete Morse function.
  Assuming the Voronoi tessellation is generated by a Poisson point process
  in $\Rspace^n$, we study the expected number of simplices in the
  $k$-dimensional weighted Delaunay mosaic as well as the expected
  number of intervals of the Morse function,
  both as functions of a radius threshold.
  As a byproduct, we obtain a new proof for the expected number of
  connected components (\emph{clumps}) in a line section of
  a circular Boolean model in $\Rspace^n$.
\end{abstract}

%\newpage

%%%%%%%%%%%%%%%%%%%%%%%%%%%%%%%%%%%%%%%%%%%%%%%%%%%%%%%%%%%%%%%%%%%%%%%%%%
%%%%%%%%%%%%%%%%%%%%%%%%%%%%%%%%%%%%%%%%%%%%%%%%%%%%%%%%%%%%%%%%%%%%%%%%%%
\section{Introduction}
\label{sec:1}
%%%%%%%%%%%%%%%%%%%%%%%%%%%%%%%%%%%%%%%%%%%%%%%%%%%%%%%%%%%%%%%%%%%%%%%%%%
%%%%%%%%%%%%%%%%%%%%%%%%%%%%%%%%%%%%%%%%%%%%%%%%%%%%%%%%%%%%%%%%%%%%%%%%%%

\begin{figure}[t]
  \centering \vspace{0.1in}
  \resizebox{!}{2.2in}{\includegraphics{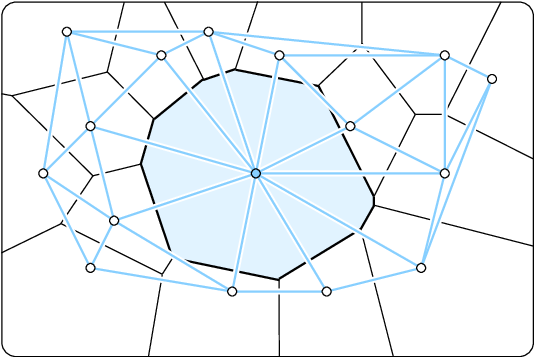}}
  \caption{Weighted Voronoi tessellation in $\Rspace^2$ with superimposed
    weighted Delaunay mosaic.
    All points have zero weight except the point with the shaded domain,
    which has positive weight.}
  \label{fig:Voronoi}
\end{figure}

Given a discrete set of points $Y \subseteq \Rspace^k$, the \emph{Voronoi tessellation}
tiles the $k$-dimensional Euclidean space with convex polyhedra,
each consisting of all points $a \in \Rspace^k$ for which a particular point
$y$ is closest among all points in $Y$.
To generalize, suppose each $y \in Y$ has a weight $w_y \in \Rspace$,
and substitute the \emph{power distance} of $a$ from $y$,
defined as $\Edist{a}{y}^2 - w_y$, for the squared Euclidean distance in
the definition of the Voronoi tessellation.
The resulting tiling of $\Rspace^k$ into convex polyhedra is known by
several names, including \emph{power diagrams} \cite{Aur87}
and \emph{Laguerre tessellations} \cite{LaZu08},
but to streamline language we will call them
\emph{weighted Voronoi tessellations}.
They do indeed look like unweighted Voronoi tessellations,
except that the hyperplane separating two neighboring polyhedra
does not necessarily lie halfway between the generating points;
see Figure \ref{fig:Voronoi}.
Our motivation for studying weighted Voronoi tessellations derives from
the extra degree of freedom --- the weight --- which permits better approximations
of observed tilings, such as cell cultures in plants \cite{PrMe16}
and microstructures of materials \cite{BWV15}.
Beyond this practical consideration,
there is an intriguing connection between the volumes of
skeleta of unweighted Voronoi tessellations and the number of simplices
in weighted Delaunay mosaics through the Crofton formula,
which is worth exploring. We will discuss it at the end of Section \ref{sec:5}.

Our preferred construction takes a $k$-dimensional slice through a
Voronoi tessellation in $\Rspace^n$; see \cite{AuIm88,Sib80}.
Specifically, if $X$ is a discrete set of points in $\Rspace^n$
and $\Rspace^k \hookrightarrow \Rspace^n$ is spanned by the first $k \leq n$
coordinate axes, then the Voronoi tessellation of $X$ in $\Rspace^n$
intersects $\Rspace^k$ in a $k$-dimensional weighted Voronoi tessellation.
The points in $\Rspace^k$ that generate the weighted tessellation
are the orthogonal projections $y_x$ of the points $x \in X$,
and their weights are $w_x = - \Edist{x}{y_x}^2$.
While all weights in this construction are non-positive,
this is not a restriction of generality because the tessellation remains
unchanged when all weights are increased by the same amount.
Indeed, every weighted Voronoi tessellation with bounded weights
can be obtained as a slice of an unweighted Voronoi tessellation.
It is often more convenient to consider the dual of a weighted Voronoi
tessellation, which is again known by several names, including
\emph{Laguerre triangulation} \cite{OBS00} and
\emph{regular triangulation} \cite{GKZ94},
but we will call them \emph{weighted Delaunay mosaics}.
An important difference to the unweighted concept is that the Voronoi
polyhedron of a weighted point may be empty, in which case this weighted
point will not be a vertex of the weighted Delaunay mosaic.
For generic sets of weighted points, the weighted Delaunay mosaic
is a simplicial complex in $\Rspace^k$.
Since we focus on slices of unweighted Voronoi tessellations,
we define the general position only in this case.
Specifically, we say a discrete set $X \subseteq \Rspace^n$ is \emph{generic}
if the following conditions are satisfied for every $0 \leq j < n$:
\smallskip \begin{enumerate}\denselist
  \item[1.] no $j+2$ points belong to a common $j$-plane,
  \item[2.] no $j+3$ points belong to a common $j$-sphere,
  \item[3.] considering the unique $j$-sphere that passes through
    $j+2$ points, no $j+1$ of these points belong to a $j$-plane
    that passes through the center of the $j$-sphere,
  \item[4.] considering the unique $j$-plane that passes through
    $j+1$ points, this plane is neither orthogonal nor parallel to $\Rspace^k$,
  \item[5.] no two points have identical distance to $\Rspace^k$.
\end{enumerate} \smallskip
For $j=0$, property 4 means that no point of $X$ is in $\Rspace^k$.
We note that the Poisson point process is generic with probability $1$.

Continuing the work started in \cite{ENR16}, we are interested
in the stochastic properties of the weighted Delaunay mosaics
and their radius functions.
To explain the latter concept, we assume the generic case in which
the mosaic is a simplicial complex, and
for every simplex $Q' \in \Delaunay{Y}$ with preimage
$Q \subseteq \Rspace^n$,
we find the smallest $(n-1)$-sphere that satisfies the following properties:
\smallskip \begin{itemize}
  \item it passes through all vertices of $Q$
    (it is a \emph{circumscribed sphere} of $Q$),
  \item the open ball it bounds does not contain any points of $X$
    (it is \emph{empty}),
  \item its center lies in $\Rspace^k$
    (it is \emph{anchored}).
\end{itemize} \smallskip
The existence of such spheres for the simplices of the weighted mosaic
can be shown in a way similar to the unweighted case \cite{Ede01}
and is left to the reader.
We call this sphere the \emph{weighted Delaunay sphere}
and its radius the \emph{weighted Delaunay radius} of $Q' \in \Delaunay{Y}$.
Similarly, when considering $Q$ instead of $Q'$, we call this sphere the \emph{anchored Delaunay sphere}
and its center the \emph{anchor} of $Q$.
The \emph{radius function} of the weighted Delaunay mosaic,
$\Rfun \colon \Delaunay{Y} \to \Rspace$, maps every simplex to its weighted Delaunay radius.
As in the unweighted case, it partitions $\Delaunay{Y}$ into \emph{intervals} of simplices
that share the same weighted Delaunay sphere and therefore the same function value \cite{BaEd15}.
These intervals have topological significance \cite{For98}:
adding the simplices in the order of increasing radius, the homotopy type
of the complex changes whenever the interval contains a single simplex
and it remains unchanged whenever the interval contains two or more simplices.
Indeed, the operation in the latter case is known as anticollapse
and has been studied extensively in combinatorial topology.
Each interval is defined by two simplices $L \subseteq U$
in the weighted Delaunay mosaic and consists of all simplices
that contain $L$ and are contained in $U$.
We call $Q' \in \Delaunay{Y}$ a \emph{critical simplex} of $\Rfun$
if it is the sole simplex in its interval: $L = Q' = U$,
and we call $Q'$ a \emph{regular simplex} of $\Rfun$, otherwise.
The \emph{type} of the interval is the pair of dimensions of the
lower and the upper bound:  $(\ell, m)$
in which $\ell = \dime{L}$ and $m = \dime{U}$.
Our main result is an extension of the stochastic findings
about the radius function of the Poisson--Delaunay mosaic in \cite{ENR16}
from the unweighted to the weighted case.
\begin{theorem}[Main Result]
  \label{thm:MainResult}
  Let $X$ be a Poisson point process with density $\density$ in $\Rspace^n$
  and $\Rspace^k \hookrightarrow \Rspace^n$.
  There are constants $\Ccon{\ell}{m}{k}{n}$ such that
  for any $r_0 \geq 0$, the expected number of intervals of type $(\ell,m)$
  in the $k$-dimensional weighted Poisson--Delaunay mosaic with center in
  a Borel set $\Region \subseteq \Rspace^k$ and weighted Delaunay radius at most $r_0$ is
  \begin{align}
    \Expected{ \ccon{\ell}{m}{k}{n} (r_0)}  &=  \Ccon{\ell}{m}{k}{n}
      \cdot \frac{\iGama{m+1-\frac{k}{n}}{\density \nu_n r_0^n}}
                 { \Gama{m+1-\frac{k}{n}}}
      \cdot \density^{\frac{k}{n}} \norm{\Region} ,
  \end{align}
  in which $\nu_n$ is the volume of the unit ball in $\Rspace^n$,
  and we give explicit computations of the constants in $k \leq 2$ dimensions.
  Similarly, the expected number of $j$-dimensional simplices in the
  weighted Poisson--Delaunay mosaic with center in
  a Borel set $\Region \subseteq \Rspace^k$ and weighted Delaunay radius at most $r_0$ is:
  \begin{align}
    \Expected{ \dcon{j}{k}{n} (r_0)}  &=  
      \left[ \sum_{m=j}^k \frac{ \iGama{m+1-\frac{k}{n}}{\density \nu_n r_0^n} }
                               {  \Gama{m+1-\frac{k}{n}} }
             \sum_{\ell=0}^j \binom{m-\ell}{m-j} \Ccon{\ell}{m}{k}{n} \right]
      \cdot  \density^{\frac{k}{n}} \norm{\Region} .
  \end{align}
\end{theorem}
Some of the values for constants $\Ccon{\ell}{m}{k}{n}$ are listed in Tables \ref{tbl:Constants1D} and \ref{tbl:Constants2D}.
In an equivalent formulation, this theorem states that the weighted Delaunay radius
of a \emph{typical interval} is Gamma-distributed,
whereas the weighted Delaunay radius of a \emph{typical simplex} is a mixture
of Gamma distributions; compare with \cite{ENR16}.
In a more general context, the contributions of this paper are to the field
of stochastic geometry, which was summarized in the
text by Schneider and Weil \cite{ScWe08}.
The particular questions on Poisson--Delaunay mosaics studied in
this paper have been pioneered by Miles almost $50$ years ago
\cite{Mil70,Mil71}.
Formulas for the weighted case have also been derived
by M{\o}ller \cite{Mol89}, but these are restricted to top-dimensional
simplices whose expected numbers can be derived using Crofton formula
and expected volumes of Voronoi skeleta.

\ourparagraph{Outline.}
Section \ref{sec:2} discusses the case $k = 1$ as a warm-up exercise.
It is sufficiently elementary so that explicit formulas can be
derived without reliance on more difficult to prove
general integral formulas.
Section \ref{sec:3} shows how to get the expected number
of connected components in the intersection of a line with
a circular Boolean model in $\Rspace^n$ using discrete Morse theory.
Section \ref{sec:4} proves a Blaschke--Petkantschin type formula
for the general weighted case.
Section \ref{sec:5} uses this formula to prove our main result.
Section \ref{sec:6} develops explicit expressions for all types
of intervals in two dimensions.
Section \ref{sec:7} concludes this paper.
Appendix \ref{app:A} introduces the special functions and distributions
used in the derivation of our results.

%%%%%%%%%%%%%%%%%%%%%%%%%%%%%%%%%%%%%%%%%%%%%%%%%%%%%%%%%%%%%%%%%%%%%%%%%%
%%%%%%%%%%%%%%%%%%%%%%%%%%%%%%%%%%%%%%%%%%%%%%%%%%%%%%%%%%%%%%%%%%%%%%%%%%
\section{One Dimension}
\label{sec:2}
%%%%%%%%%%%%%%%%%%%%%%%%%%%%%%%%%%%%%%%%%%%%%%%%%%%%%%%%%%%%%%%%%%%%%%%%%%
%%%%%%%%%%%%%%%%%%%%%%%%%%%%%%%%%%%%%%%%%%%%%%%%%%%%%%%%%%%%%%%%%%%%%%%%%%

In $k = 1$ dimension, the weighted Delaunay mosaic has a simple
structure so that results can be obtained by elementary means.

\begin{figure}[t]
  \centering \vspace{0.1in}
    \resizebox{0.48\textwidth}{!}{\includegraphics{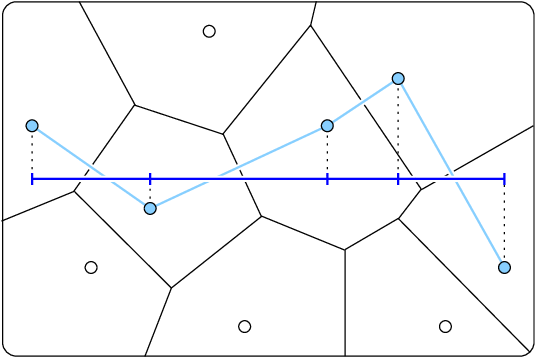}} \hspace{0.1in}
    \resizebox{0.48\textwidth}{!}{\includegraphics{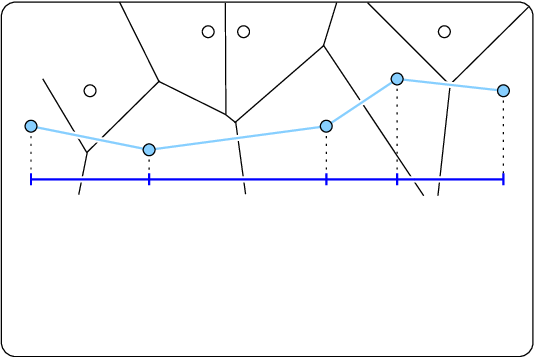}}
  \caption{\emph{Left:} a $1$-dimensional weighted Voronoi tessellation
    as a slice of a $2$-dimensional unweighted Voronoi tessellation.
    The weighted Delaunay mosaic in $\Rspace^1$ is the projection
    of a chain of edges in the $2$-dimensional unweighted Delaunay mosaic.
    \emph{Right:} reflecting the points across $\Rspace^1$ affects
    the $2$-dimensional Voronoi tessellation but not the $1$-dimensional slice.}
  \label{fig:projection}
\end{figure}

\ourparagraph{Slice construction.}
Let $n \geq 2$ and let $X \subseteq \Rspace^n$ be a stationary
Poisson point process with density $\density > 0$.
We write $\Rspace^1 \hookrightarrow \Rspace^n$ for the first
coordinate axis, which is a directed line passing through $\Rspace^n$.
For each point $x = (x_1, x_2, \ldots, x_n) \in X$,
we write $y_x = (x_1, 0, \ldots, 0)$ for the projection onto $\Rspace^1$
and $- w_x = x_2^2 + x_3^2 + \ldots + x_n^2$ for its squared distance
from the line.
Letting $Y = \{ (y_x, w_x) \mid x \in X \}$
be the resulting set of weighted points in $\Rspace^1$,
we are interested in its weighted Voronoi tessellation, $\Voronoi{Y}$,
and its weighted Delaunay mosaic, $\Delaunay{Y}$.
By construction, the former is the intersection of the $n$-dimensional
(unweighted) Voronoi tessellation with the line:
$\Voronoi{Y} = \{ \voronoi{x} \cap \Rspace^1 \mid x \in X \}$.
As discussed above, the interval $\voronoi{x} \cap \Rspace^1$ belongs
to the weighted Voronoi tessellation iff there is an anchored Delaunay sphere of $x$,
that is: an empty sphere centered in $\Rspace^1$ that passes through $x$.
Similarly, two weighted Voronoi domains, $\voronoi{x} \cap \Rspace^1$
and $\voronoi{u} \cap \Rspace^1$, share an endpoint iff 
there is an empty anchored Delaunay sphere passing through $x$ and $u$.
It follows that every edge in $\Delaunay{Y}$ is the projection
of an edge in $\Delaunay{X}$;
see Figure \ref{fig:projection}.
As suggested in this figure, we can simplify the
construction by reducing $n$ to $2$.
Writing $\UpperHalfPlane$ for the half-plane of points whose first coordinate
is arbitrary, whose second coordinate is non-negative,
and whose remaining $n-2$ coordinates vanish,
we map $x = (x_1, x_2, \ldots, x_n) \in \Rspace^n$
to $x' = (x_1, \sum_{i=2}^n x_i^2, 0, \ldots, 0) \in \UpperHalfPlane$.
This amounts to rotating $x$ about $\Rspace^1$ into $\UpperHalfPlane$.
Let $X'$ be the resulting set of points in $\UpperHalfPlane$
and $Y'$ the set of weighted points in $\Rspace^1$
obtained by projection from $X'$.
Then $Y = Y'$, which shows that $X$ and $X'$ define the same
$1$-dimensional weighted Voronoi tessellation and weighted Delaunay mosaic.
There is a small price to pay for the simplification,
namely that the projected Poisson point process in $\UpperHalfPlane$ is
not necessarily homogeneous.
Specifically, the projected process in $\UpperHalfPlane$
is a Poisson point process with \emph{intensity}
$\intensity{x} = \sigma_{n-1} \density x_2^{n-2}$,
in which $\sigma_{n-1}$ is the $(n-2)$-dimensional volume
of the unit sphere in $\Rspace^{n-1}$.

\ourparagraph{Interval structure.}
We now return to the intervals of the radius function in
one dimension, $\Rfun \colon \Delaunay{Y} \to \Rspace$.
In the assumed generic case, $\Delaunay{Y}$ contains only two kinds
of simplices:  vertices and edges.
By definition, the value of $\Rfun$ at a simplex
$Q' \in \Delaunay{Y}$ is the radius of the anchored
Delaunay sphere of the preimage of $Q'$.
There are only three types of intervals $[L, U]$:
\smallskip \begin{description}\denselist
  \item[$(0,0)$:] here $L = U$ and $\dime{L} = \dime{U} = 0$.
    The interval contains a single and therefore critical vertex.
  \item[$(1,1)$:] here $L = U$ and $\dime{L} = \dime{U} = 1$.
    The interval contains a single and therefore critical edge.
  \item[$(0,1)$:] here $L \subseteq U$ and $\dime{L} = \dime{U} - 1$.
    The interval is a pair consisting of a regular vertex and a regular edge.
    We call it a \emph{vertex-edge pair} if the vertex precedes the edge
    as we go from left to right, and we call it an \emph{edge-vertex pair},
    otherwise.
\end{description} \smallskip
\begin{figure}[t]
  \centering \vspace{0.1in}
    \resizebox{!}{1.8in}{\includegraphics{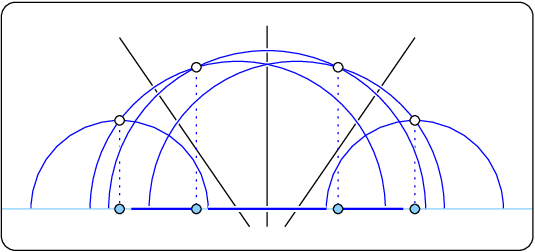}} \hspace{0.2in}
  \caption{From \emph{left} to \emph{right} on the horizontal line:
    a critical vertex, an edge-vertex pair, a critical edge,
    a vertex-edge pair, and another critical vertex.}
  \label{fig:1D-projection}
\end{figure}
The cases can be distinguished geometrically, as illustrated
in Figure \ref{fig:1D-projection}.
Let $x = (x_1, x_2) \in \UpperHalfPlane$ and $y_x = (x_1, 0)$ with weight
$w_x = - x_2^2$.
Then $L = U= \{y_x\}$ is a critical vertex of $\Delaunay{Y}$
iff $y_x$ is the anchor of $x$.
Otherwise, the anchored Delaunay circle of $x$ also passes through a second point,
$u \in X \subseteq \UpperHalfPlane$, with $y_x$ and $y_u$ on the same side
of the anchor.
In this case, $L = \{y_x\}$ and $U = \{y_x, y_u\}$ form a vertex-edge
or an edge-vertex pair.
Finally, we have a critical edge $L = U = \{y_x, y_u\}$
if $y_x$ and $y_u$ lie on opposite sides of the anchor.

We will make essential use of the geometric characterization of
interval types when we compute their expected numbers.
To simplify the computation, we note that the structure
along $\Rspace^1$ is a strict repetition of the following pattern:
a critical vertex, a non-negative number of edge-vertex pairs,
a critical edge, and a non-negative number of vertex-edge pairs.

\ourparagraph{Critical vertices.}
We begin with computing the number of critical vertices,
$\ccon{0}{0}{1}{n}$, inside a region $\Region \subseteq \Rspace^1$
and with weighted Delaunay radius at most some threshold $r_0$.
Let $x = (x_1, x_2) \in X \subseteq \UpperHalfPlane$ and note that the
smallest anchored circle passing through $x$ has center
$y_x = (x_1, 0)$ and radius $r = x_2$.
Write $\Pempty{x}$ for the probability that this circle is empty,
$\One_{\Region}(x)$ for the indicator that $y_x \in \Region$,
and $\One_{r_0}(x)$ for the indicator that $r \leq r_0$.
We use the Slivnyak--Mecke formula to compute
\begin{align}
  \Expected{\ccon{0}{0}{1}{n} (r_0)} 
    &=  \int\displaylimits_{x \in \UpperHalfPlane} \One_\Region(x) \One_{r_0}(x)
          \Pempty{x} \intensity{x} \diff x ;
\end{align}
compare with \cite{ENR16}.
The intensity measure of the upper semi-circle with radius $r$ is
of course $\density$ times the volume of an $n$-ball with radius $r$,
which we write as $\density \nu_r r^n$.
Hence, $\Pempty{x} = e^{- \density \nu_n r^n}$.
In other words, the probability that the anchored
circle is empty is the probability that
the $n$-ball whose points get rotated into the semi-disk
is empty.
So we have
\begin{align}
  \Expected{\ccon{0}{0}{1}{n} (r_0)} 
    &=  \int\displaylimits_{x_1 \in \Region} \int\displaylimits_{r=0}^{r_0}
          e^{- \density \nu_n r^n} \density \sigma_{n-1} r^{n-2}
          \diff r \diff x_1
     =  \norm{\Region} \sigma_{n-1} \density
        \int\displaylimits_{r=0}^{r_0} r^{n-2} e^{- \density \nu_n r^n} \diff r.
  \label{eqn:criticalvertices1D}
\end{align}
To evaluate this integral, we use the identity on Gamma functions
proved as Lemma \ref{lem:GammaFunction} in Appendix \ref{app:A},
where the functions are defined.
In this application,
the integral on the right-hand side in \eqref{eqn:criticalvertices1D}
evaluates to $\iGama{1-\frac{1}{n}}{\density \nu_n r_0^n}
            / [n \cdot (\density \nu_n)^{1-\frac{1}{n}}]$.
Writing $\ccon{0}{0}{1}{n} = \ccon{0}{0}{1}{n} (\infty)$,
we set $r_0 = \infty$ to get the expected total number of
critical vertices, and we write the expected number up to
weighted Delaunay radius $r_0$ as a fraction of the former:
\begin{align}
  \Expected{\ccon{0}{0}{1}{n}}
    &= \frac{\sigma_{n-1} \Gama{1-\frac{1}{n}}}{n \nu_n^{1-1/n}}
       \cdot \norm{\Region} \density^{\frac{1}{n}} ,
  \label{eqn:critical-vertices-1} \\
  \Expected{\ccon{0}{0}{1}{n} (r_0)}  &= 
    \frac{\iGama{1-\frac{1}{n}}{\density \nu_n r_0^n}}
         {\Gama{1-\frac{1}{n}}}
    \cdot \Expected{\ccon{0}{0}{1}{n}} .
  \label{eqn:critical-vertices-2}
\end{align}

\ourparagraph{Regular edges.}
To count the regular edges --- or intervals of type $(0,1)$ ---
we again use the Slivnyak--Mecke formula.
Let $x = (x_1, x_2)$ and $u = (u_1, u_2)$ be two points
in $X \subseteq \UpperHalfPlane$.
There is a unique anchored circle that passes through both points,
and the edge connecting $y_x$ and $y_u$ belongs to $\Delaunay{Y}$
iff this circle is empty.
Writing $(z_1, 0)$ for the center and $r$ for the radius,
the edge is critical, if $x_1 < z_1 < u_1$,
and regular, otherwise; see Figure \ref{fig:1D-projection}.
Write $\Pempty{x,u}$ for the probability that the unique
anchored circle passing through $x$ and $u$ is empty,
$\One_{\Region} (x,u)$ for the indicator that $z_1 \in \Region$,
write $\One_{r_0}(x,u)$ for the indicator that $r \leq r_0$, and
$\One_{0,1}(x,u)$ for the indicator that
$x_1$ and $u_1$ lie on the same side of $z_1$.
By Slivnyak--Mecke formula, we have
\begin{align}
  \Expected{\ccon{0}{1}{1}{n} (r_0)}  &=  \frac{1}{2!}
    \int\displaylimits_{u \in \UpperHalfPlane}
    \int\displaylimits_{x \in \UpperHalfPlane}
		\One_\Region (x,u) \One_{r_0} (x,u)
      \One_{0,1} (x,u) \Pempty{x,u} \intensity{x} \intensity{u}
      \diff x \diff u .
  \label{eqn:regular-edge}
\end{align}
We already know that $\Pempty{x,u} = e^{- \density \nu_n r^n}$.
To compute the rest, we do a change of variables,
re-parametrizing the points by the center and radius of the unique anchored
circle passing through them and two angles:
$x  =  (z_1 + r \cos \xi, r \sin \xi)$
and $u  =  (z_1 + r \cos \upsilon, r \sin \upsilon)$,
in which $0 \leq \xi, \upsilon < \pi$.
This is a bijection up to a set of measure $0$.
The Jacobian of this change of variables is the absolute determinant
of the matrix of old variables derived by the new variables:
\begin{align}
  J  &=  \abs{\left| \begin{array}{cccc}
                1  &  \cos \xi       & - r \sin \xi  &  0                 \\
                0  &  \sin \xi       &   r \cos \xi  &  0                 \\
                1  &  \cos \upsilon  &   0           &  - r \sin \upsilon \\
                0  &  \sin \upsilon  &   0           &    r \cos \upsilon
         \end{array} \right|}
      =  r^2 | \cos \upsilon - \cos \xi | .
\end{align}
With the new variables, the indicators can be absorbed into
integration limits:
$\One_\Region (x,u) = 1$ iff $z_1 \in \Region$,
and $\One_{0,1} (x,u) = 1$ iff $\xi$ and $\upsilon$ are either
both smaller or both larger than $\frac{\pi}{2}$.
The two cases are symmetric, so we assume the former and
multiply with $2$.
The integral in \eqref{eqn:regular-edge} thus turns into
\begin{align}
  \Expected{\ccon{0}{1}{1}{n} (r_0)}
    &=  \!\!\int\displaylimits_{z_1 \in \Region}
        \int\displaylimits_{r=0}^{r_0}\!\! e^{- \density \nu_n r^n}
        \!\!\!\!\!\!\!\int\displaylimits_{0 \leq \xi,
                                      \upsilon < \frac{\pi}{2}}\!\!\!\!\!\!\!
          \density^2 \sigma_{n-1}^2 (r^2 \sin \xi \sin \upsilon)^{n-2}
          r^2 |\cos \upsilon - \cos \xi| \diff \xi \diff \upsilon \diff r \diff z_1
  \label{eqn:regular-edges-1} \\
    &=  \norm{\Region} \density^2 \sigma_{n-1}^2
        \!\int\displaylimits_{r=0}^{r_0}\!
          e^{- \density \nu_n r^n} r^{2n-2} \diff r
        \!\!\!\!\!\!\!\int\displaylimits_{0 \leq \xi,
                                          \upsilon < \frac{\pi}{2}}\!\!\!\!\!\!\!
          (\sin \xi \sin \upsilon)^{n-2}
          |\cos \upsilon - \cos \xi| \diff \xi \diff \upsilon .
  \label{eqn:regular-edges-2}
\end{align}
We apply Lemma \ref{lem:GammaFunction} to evaluate the integral
over the radius, and we use the Mathematica software to evaluate
the integral over the two angles:
\begin{align}
  \int\displaylimits_{r \leq r_0} r^{2n-2} e^{- \density \nu_n r^n} \diff r
    &=  \tfrac{ \iGama{2-\frac{1}{n}}{\density \nu_n r_0^n} }
             { n (\density \nu_n)^{2 - \frac{1}{n}} } ,          \\
  \int\displaylimits_{0 \leq \xi, \upsilon < \frac{\pi}{2}}
        (\sin \xi \sin \upsilon)^{n-2} |\cos \upsilon - \cos \xi|
        \diff \xi \diff \upsilon
    &=  \tfrac{\sqrt{\pi}}{n-1}
        \left[ \tfrac{2 \Gama{n-1}}{\Gama{n-\frac{1}{2}}}
             - \tfrac{\Gama{\frac{n-1}{2}}}{\Gama{\frac{n}{2}}} \right].
\end{align}
Setting $r_0 = \infty$, we get the expected total number of regular edges,
and as before we write the expected number up to weighted Delaunay radius $r_0$ as a
fraction of the total number:
\begin{align}
  \Expected{\ccon{0}{1}{1}{n}} 
    &=  \tfrac{\sigma_{n-1}^2 \Gama{2-\frac{1}{n}}}{n \nu_n^{2-1/n}}
        \tfrac{\sqrt{\pi}}{n-1}
        \left[ \tfrac{2 \Gama{n-1}}{\Gama{n-\frac{1}{2}}}
             - \tfrac{\Gama{\frac{n-1}{2}}}{\Gama{\frac{n}{2}}} \right]
        \cdot \norm{\Region} \density^{\frac{1}{n}} , 
  \label{eqn:regular-edges-3} \\
  \Expected{\ccon{0}{1}{1}{n} (r_0)}
    &=   \tfrac{\iGama{2-\frac{1}{n}}{\density \nu_n r_0^n}}
              {\Gama{2-\frac{1}{n}}}
         \cdot \Expected{\ccon{0}{1}{1}{n}}.
  \label{eqn:regular-edges-4}
\end{align}

\ourparagraph{Summary.}
Recall that the critical vertices and the critical edges alternate
along $\Rspace^1$,
which implies that their expected total number is the same.
The dependence on the radius threshold, $r_0$, is however different.
Here we notice that the dependence on the radius for
$\ccon{1}{1}{1}{n}$ is the same as for $\ccon{0}{1}{1}{n}$
because what changes in the integration are only the
admissible angles.
Extracting the constants from the formulas for the expectation,
we use \eqref{eqn:critical-vertices-1} and \eqref{eqn:regular-edges-3} to get
\begin{align}
  \Ccon{0}{0}{1}{n}   =  \Ccon{1}{1}{1}{n}
                     &=  \tfrac{\sigma_{n-1} \Gama{1-\frac{1}{n}}}{n \nu_n^{1-1/n}} ,
    \label{eqn:constant00} \\
  \Ccon{0}{1}{1}{n}  &=  \tfrac{\sigma_{n-1}^2 \sqrt{\pi} \Gama{2-\frac{1}{n}}}
                              {n (n-1) \nu_n^{2-1/n}}
                         \left[ \tfrac{2 \Gama{n-1}}{\Gama{n-\frac{1}{2}}}
                              - \tfrac{\Gama{\frac{n-1}{2}}}{\Gama{\frac{n}{2}}}
                                                                            \right] ;
    \label{eqn:constant01}
\end{align}
see Table \ref{tbl:Constants1D}.
We write the expectations as fractions of these constants times
the size of the region times the $n$-th root of the density in $\Rspace^n$:
\begin{align}
  \Expected{\ccon{0}{0}{1}{n} (r_0)}  &=  \Ccon{0}{0}{1}{n} \cdot
    \tfrac{\iGama{1-\frac{1}{n}}{\density \nu_n r_0^n}}
          {\Gama{1-\frac{1}{n}}} \cdot \norm{\Region} \density^{1/n} ,      \\
  \Expected{\ccon{0}{1}{1}{n} (r_0)}  &=  \Ccon{0}{1}{1}{n} \cdot
    \tfrac{\iGama{2-\frac{1}{n}}{\density \nu_n r_0^n}}
          {\Gama{2-\frac{1}{n}}} \cdot \norm{\Region} \density^{1/n} ,      \\
  \Expected{\ccon{1}{1}{1}{n} (r_0)}  &=  \Ccon{1}{1}{1}{n} \cdot
    \tfrac{\iGama{2-\frac{1}{n}}{\density \nu_n r_0^n}}
          {\Gama{2-\frac{1}{n}}} \cdot \norm{\Region} \density^{1/n} .
\end{align}
To get the corresponding results for the simplices in the weighted
Delaunay mosaic, we note that the number of vertices is
$\dcon{0}{1}{n} = \ccon{0}{0}{1}{n} + \ccon{0}{1}{1}{n}$
and the number of edges is
$\dcon{1}{1}{n} = \ccon{0}{1}{1}{n} + \ccon{1}{1}{1}{n}$.
The two are the same, but this is not true
if we limit the radius to a finite threshold.
Indeed, the radius of a typical edge is Gamma distributed
while the radius of a typical vertex follows a linear combination
of two Gamma distributions.
In the limit, when $n \to \infty$, the constants are
$\lim_{n \to \infty} \Ccon{0}{0}{1}{n} = \sqrt{e}$,
$\lim_{n \to \infty} \Ccon{0}{1}{1}{n} = \sqrt{e} (\sqrt{2} - 1)$, and
$\lim_{n \to \infty} \Dcon{0}{1}{n} = \lim_{n \to \infty} \Dcon{1}{1}{n} = \sqrt{2 e}$,
which can again be verified using the Mathematica software.
\begin{table}[hbt]
  \centering
  \small{ \begin{tabular}{r||rrrr rrrr rrrr}
    & \multicolumn{1}{c}{$n=2$}& \multicolumn{1}{c}{$3$} & \multicolumn{1}{c}{$4$}   
    & \multicolumn{1}{c}{$5$}  & \multicolumn{1}{c}{$6$} & \multicolumn{1}{c}{$7$}
    & \multicolumn{1}{c}{$8$}  & \multicolumn{1}{c}{$9$}& \multicolumn{1}{c}{$\ldots$}
    & \multicolumn{1}{c}{$20$} & \multicolumn{1}{c}{$\ldots$}& \multicolumn{1}{c}{$\infty$}
                                                                          \\ \hline \hline
    $\Ccon{0}{0}{1}{n}$
      & $1.00$ & $1.09$ & $1.16$ & $1.22$ & $1.26$ & $1.29$ & $1.32$ & $1.35$
                                 &$\ldots$& $1.47$ &$\ldots$& $1.65$  \\
    $\Ccon{0}{1}{1}{n}$
      & $0.27$ & $0.36$ & $0.42$ & $0.45$ & $0.48$ & $0.50$ & $0.51$ & $0.53$
                                 &$\ldots$& $0.60$ &$\ldots$& $0.68$  \\ \hline
    $\Dcon{0}{1}{n}$
      & $1.27$ & $1.46$ & $1.58$ & $1.67$ & $1.74$ & $1.79$ & $1.84$ & $1.87$
                                 &$\ldots$& $2.07$ &$\ldots$& $2.33$
%%  $\Ccon{0}{0}{1}{n}/\Ccon{0}{1}{1}{n}$
%%    & $3.66$ & $3.00$ & $2.79$ & $2.69$ & $2.63$ & $2.60$ & $2.57$
%%                               &$\ldots$& $2.47$ &$\ldots$& $2.41$
  \end{tabular} }
  \caption{The rounded constants in the expressions of the expected number of
    intervals and simplices of a $1$-dimensional weighted Delaunay mosaic.
    The ratio of the expected number of critical edges over
    the expected number of regular edges it is monotonically decreasing.
    It follows that we can infer the ambient dimension from the ratio.}
  \label{tbl:Constants1D}
\end{table}

%%%%%%%%%%%%%%%%%%%%%%%%%%%%%%%%%%%%%%%%%%%%%%%%%%%%%%%%%%%%%%%%%%%%%%%%%%
%%%%%%%%%%%%%%%%%%%%%%%%%%%%%%%%%%%%%%%%%%%%%%%%%%%%%%%%%%%%%%%%%%%%%%%%%%
\section{Connection to Boolean Model}
\label{sec:3}
%%%%%%%%%%%%%%%%%%%%%%%%%%%%%%%%%%%%%%%%%%%%%%%%%%%%%%%%%%%%%%%%%%%%%%%%%%
%%%%%%%%%%%%%%%%%%%%%%%%%%%%%%%%%%%%%%%%%%%%%%%%%%%%%%%%%%%%%%%%%%%%%%%%%%

\begin{figure}[hb]
  \centering \vspace{0.05in}
  \resizebox{!}{1.6in}{\includegraphics{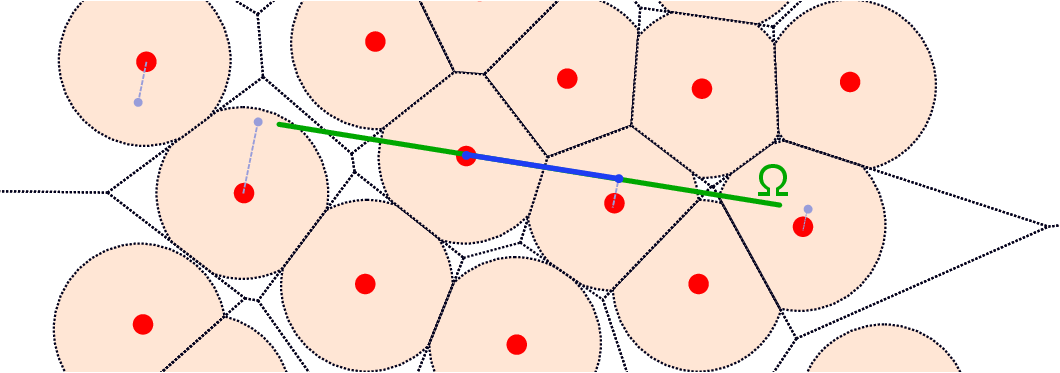}} \hspace{0.2in}
  \caption{Intersection of a union of $2$-dimensional balls with a line segment, $\Region$.
    This intersection has three components, two more than the restricted
    weighted Delaunay mosaic, which consists of two vertices and the connecting
    edge in the middle of $\Region$.
    The restricted mosaic misses the two tail components because the centers of the
    corresponding balls do not project into $\Region$.}
  \label{fig:restricted}
\end{figure}

Let $X$ be a Poisson point process with density $\density$ in $\Rspace^n$,
and write $X_r$ for the union of closed balls of fixed radius $r$
whose centers are in $X$.
This random set is sometimes referred to as the \emph{Boolean model} \cite{ScWe08}.
Let $\Region \subseteq \Rspace^1 \subseteq \Rspace^n$ be a line segment,
and consider $X_r \cap \Region$.
We are interested in the  connected components in this intersection
and claim that their number satisfies
$\beta_0(\Del{r}{Y; \Region}) \leq \beta_0 (X_r \cap \Region) \leq \beta_0(\Del{r}{Y; \Region}) + 2$,
in which $\Del{r}{Y; \Region}$ is the subcomplex of the weighted Delaunay mosaic
that consists of all simplices with radius at most $r$
whose weighted Delaunay center lies in $\Region$.
This follows from the general observation that the weighted Delaunay mosaic of
a set of points $y \in \Rspace^k$ with weights $w_y$ is homotopy equivalent to the
union of power balls,
$Y_r = \{ a \in \Rspace^k \mid \Edist{a}{y}^2 - w_y \leq r^2 \}$,
and $Y_r \cap \Region = X_r \cap \Region$.
Indeed, the weighted Delaunay complex can be defined as the nerve of
the decomposition of $Y_r$ with the weighted Voronoi tessellation,
so the Nerve Theorem asserts the homotopy equivalence;
see \cite{EdHa10} for details.
By restricting the Delaunay mosaic to a line segment,
we can lose up to two connected components
at the ends of $\Region$; see Figure \ref{fig:restricted}.

\begin{figure}[ht]
  \centering \vspace{0.1in}
  \resizebox{!}{2in}{\includegraphics{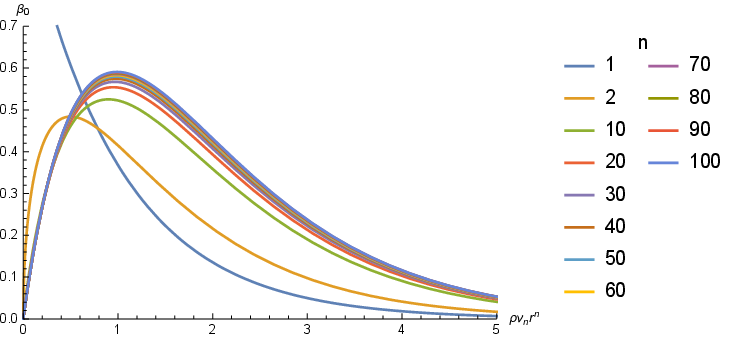}} \hspace{0.2in}
  \caption[Expected number of connected components in a line section per unit length]%
    {The expected number of connected components per unit length as a function
    of the radius.
    To facilitate the comparison of the graphs in different dimensions, $n$,
    we rescale such that a unit along the horizontal axis is %$\density \nu_n r^n$,
    the expected number of points inside a ball of radius $r$ in $\Rspace^n$.
    %After rescaling, the functions are
    %$\beta_0 (\density \nu_n r^n) = \Expected{\beta_0 (X_r \cap \Region)}$
    %with $\norm{\Region} = 1$.}
		}
  \label{fig:betti}
\end{figure}

Following the evolution of the nested complexes $\Del{r}{Y; \Region}$,
as $r$ goes from $0$ to $\infty$,
we observe that upon entering the complex
a critical vertex creates a new component,
a regular interval does not affect the homotopy type,
and a critical edge connects two components;
compare with Figure \ref{fig:1D-projection}.
It follows that the expected number of components in $\Del{r}{Y; \Region}$ is
\begin{align}
  \Expected{\ccon{0}{0}{1}{n} (r) - \ccon{1}{1}{1}{n} (r)} 
     &= \tfrac{\sigma_{n-1} \Gama{1-\frac{1}{n}}}{n \nu_n^{1-1/n}} 
        \left[ \tfrac{\iGama{1-\frac{1}{n}}{\density \nu_n r^n}}{\Gama{1-\frac{1}{n}}}
             - \tfrac{\iGama{2-\frac{1}{n}}{\density \nu_n r^n}}{\Gama{2-\frac{1}{n}}} \right]
        \cdot \density^{\frac{1}{n}} \norm{\Region} \\
     &= \tfrac{\sigma_{n-1}}{n \nu_n^{1-1/n}} 
        \left[ \iGama{1-\tfrac{1}{n}}{\density \nu_n r^n}
             - \tfrac{\iGama{2-\frac{1}{n}}{\density \nu_n r^n}}{{1-\frac{1}{n}}} \right]
        \cdot \density^{\frac{1}{n}} \norm{\Region} .
  \label{eq:boolean-1}
\end{align}
We write $A = \density \nu_n r^n$, use the definition of the incomplete Gamma function,
and integrate by parts to get
\begin{align}
  \iGama{2-\frac{1}{n}}{A} &= \int\displaylimits_0^A x^{1-\frac{1}{n}} e^{-x} \diff x 
     = \left[ - x^{1-\frac{1}{n}} e^{-x} \right]_0^A + \left( 1-\tfrac{1}{n} \right)
       \int\displaylimits_0^A x^{-\frac{1}{n}} e^{-x} \diff x \\
    &= -A^{1-\frac{1}{n}} e^{-A} + \left( 1-\tfrac{1}{n} \right) \iGama{1-\tfrac{1}{n}}{A}.
  \label{eq:bool-gamma}
\end{align}
Noticing that
$A^{1-\frac{1}{n}} \density^{1/n} = (\density \nu_n r^n)^{1-\frac{1}{n}} \density^{1/n}
                                  = \density \nu_n^{1-\frac{1}{n}} r^{n-1}$,
we plug \eqref{eq:bool-gamma} into \eqref{eq:boolean-1} to obtain
\begin{align}
  \Expected{\beta_0(\Del{r}{Y; \Region})}
    &= \tfrac{\sigma_{n-1}}{n \nu_n^{1-1/n}} \tfrac{1}{1-\frac{1}{n}}
       e^{-\density \nu_n r^n} \density \nu_n^{1-\frac{1}{n}} r^{n-1} \norm{\Region}
     = \tfrac{\sigma_{n-1}}{n-1} r^{n-1} e^{-\density \nu_n r^n} \density \norm{\Region} \\
    &= \nu_{n-1} r^{n-1} e^{-\density \nu_n r^n} \density \norm{\Region},
  \label{eq:boolean}
\end{align}
where we use the identity $\tfrac{\sigma_{n-1}}{n-1} = \nu_{n-1}$ in the last transition.
In short, \eqref{eq:boolean} gives an explicit formula for the expected density
of connected components in the Boolean model in $\Rspace^n$
intersected with a line.
While the authors did not find the explicit expression in the literature,
this result is not new and follows after some straightforward computations from
\cite[Excercise 4.8]{Hall88}.
Our aim is to provide another, more topological view on the problem.
The graphs of $\beta_0$ for different dimensions $n$ are shown in Figure \ref{fig:betti}.
Using Crofton formula \cite[Theorem 9.4.7]{ScWe08} but see also \cite{Had57}
and the fact that almost every connected component is a line segment
that meets the boundary of the Boolean model in two points,
\eqref{eq:boolean} can be transformed into a statement
about the boundary of $X_r$:
\begin{align}
  \overline{V}_{n-1}(X_r) &= 2 \sqrt{\pi} \tfrac{\Gama{\frac{n}{2}}}{\Gama{\frac{n+1}{2}}}
                             \nu_{n-1} r^{n-1} e^{-\density \nu_n r^n} \density,
\end{align}
in which $\overline{V}_{n-1}(X_r)$ is the expected density of
$(n-1)$-dimensional volume of the boundary;
see \cite[Section 9]{ScWe08} for the detailed discussion of the quantity.

%%%%%%%%%%%%%%%%%%%%%%%%%%%%%%%%%%%%%%%%%%%%%%%%%%%%%%%%%%%%%%%%%%%%%%%%%%
%%%%%%%%%%%%%%%%%%%%%%%%%%%%%%%%%%%%%%%%%%%%%%%%%%%%%%%%%%%%%%%%%%%%%%%%%%
\section{Anchored Blaschke--Petkantschin Formula}
\label{sec:4}
%%%%%%%%%%%%%%%%%%%%%%%%%%%%%%%%%%%%%%%%%%%%%%%%%%%%%%%%%%%%%%%%%%%%%%%%%%
%%%%%%%%%%%%%%%%%%%%%%%%%%%%%%%%%%%%%%%%%%%%%%%%%%%%%%%%%%%%%%%%%%%%%%%%%%

To extend the results in the previous section from $1$ to $k$ dimensions,
we first generalize the Blaschke--Petkantschin formula for spheres
stated as Theorem 7.3.1 in \cite{ScWe08}.

\ourparagraph{Setting the stage.}
Recall that $k \leq n$ are positive integers, and that we write
$\Rspace^k$ for the $k$-dimensional linear subspace spanned by
the first $k$ coordinate vectors of $\Rspace^n$.
While we used uppercase letters to denote simplices in the previous sections,
we now write $\xxx$ for a sequence of $m+1 \leq k+1$ points in $\Rspace^n$.
The reason for the change of notation is that we integrate over
all such sequences and do not limit ourselves to points in the
Poisson point process.
Similarly, we write $\uuu$ if the $m+1$ points lie on
the unit sphere.
As usual, we do not distinguish between a simplex and its vertices,
so we write $\Vol{m}{\xxx}$ for the $m$-dimensional Lebesgue measure
of the convex hull of $\xxx$.
Assuming the $m+1$ points are in general position in $\Rspace^n$,
the affine hull of $\xxx$ is an $m$-plane, $M = \aff{\xxx}$.
Furthermore, the set of centers of the spheres that pass through all points of $\xxx$
is an $(n-m)$-plane, $M^\perp$, orthogonal to $M$.
Generically, the intersection of $M^\perp$ with $\Rspace^k$ is a plane
of dimension $k-m$.
The center of the smallest anchored sphere passing through $\xxx$
is the point of this intersection that is the closest to $\xxx$.

\ourparagraph{Top-dimensional case.}
We first show how to transform an integral over $m+1 = k+1$ points
to the integral over the unique anchored sphere passing through these points.
\begin{lemma}[Blaschke--Petkantschin for Top-dimensional Simplices]
  \label{lem:BPforTopdimensionalSimplices}
  Let $0 \leq k \leq n$.
  Then every measurable non-negative function
  $f \colon (\Rspace^n)^{k+1} \to \Rspace$ satisfies
  \begin{align}
    \int\displaylimits_{\xxx \in (\Rspace^n)^{k+1}} f(\xxx) \diff \xxx
      &=  \int\displaylimits_{y \in \Rspace^k}
          \int\displaylimits_{r \geq 0}
          \int\displaylimits_{\uuu \in (\Sspace^{n-1})^{k+1}}
            f(y+r\uuu) r^{(n-1)(k+1)} k! \Vol{k}{\uuu'}
            \diff \uuu \diff r \diff y ,
    \label{eqn:BP1}
  \end{align}
  in which $\uuu'$ is the projection of $\uuu$ to $\Rspace^k$,
  $\Vol{k}{\uuu'}$ is the Lebesgue measure of the $k$-simplex,
  and we use the standard spherical measure on $\Sspace^{n-1}$.
\end{lemma}
\ourproof
  We follow the proof of Theorem 7.3.1 in \cite{ScWe08},
  with just slight modifications.
  Recall first that we choose the coordinates in $\Rspace^n$
  so that the projection of $x = (x_1, x_2, \ldots, x_n)$ to
  $\Rspace^k \hookrightarrow \Rspace^n$
  is $x' = (x_1, \ldots, x_k, 0, \ldots, 0)$.
  The claimed relation is a change of variables:
  on the right-hand side, we represent the points $\xxx$ by the
  center $y  \in \Rspace^k\hookrightarrow \Rspace^n$ of the anchored sphere
  passing through these points,
  its radius $r$,
  and $k$ points $\uuu$ on the unit sphere $\Sspace^{n-1} \hookrightarrow \Rspace^n$.
  This change of variables is the mapping
  $\varphi \colon \Rspace^k \times [0,\infty) \times (\Sspace^{n-1})^{k+1}
                  \to (\Rspace^n)^{k+1}$
  defined by
  $\varphi (y, r, \uuu_0, \uuu_1, \ldots, \uuu_k)
    = (y+r \uuu_0, y+r \uuu_1, \ldots, y+r \uuu_k)$,
  we note that $\varphi$ is bijective up to a measure $0$ subset of the domain.
  We claim the Jacobian of $\varphi$ is
  \begin{align}
    J(y, r, \uuu)  &=  r^{(n-1)(k+1)} k! \Vol{k}{\uuu'} ,
    \label{eqn:Jacobian}
  \end{align}
  in which $\uuu' = (\uuu_0', \uuu'_1, \ldots, \uuu'_k)$
  is the projection of $\uuu$ to $\Rspace^k$.
  To prove \eqref{eqn:Jacobian} at a particular point $(y, r, \uuu)$,
  we choose local coordinates around every point $\uuu_i$ on the sphere.
  We choose them such that the matrix $[\uuu_i \dot{\uuu}_i]$ is orthogonal,
  for every $0 \leq i \leq k$,
  in which $\dot{\uuu}_i$ is the $n \times (n-1)$ matrix of partial
  derivatives with respect to the $n-1$ local coordinates.
  This is the same parametrization as in \cite{ScWe08}.
  With this, the Jacobian is the absolute value of the
  $n(k+1) \times n(k+1)$ determinant:
  \begin{align}
    J (y, r, \uuu)
%       &=  \left|
%      \begin{array}{cc:cccc}
%        E_{k,k} & \uuu_0'      & r \dot{\uuu}_0'      & 0 & \ldots & 0 \\
%        0       & \uuu_0^\perp & r \dot{\uuu}_0^\perp & 0 & \ldots & 0 \\
%        E_{k,k} & \uuu_1'      & 0 & r \dot{\uuu}_1'      & \ldots & 0 \\
%        0       & \uuu_1^\perp & 0 & r \dot{\uuu}_1^\perp & \ldots & 0 \\
%        \vdots     & \vdots & \vdots & \vdots             & \ddots & \vdots \\
%        E_{k,k} & \uuu_k'      & 0 & 0 & \ldots & r \dot{\uuu}_k'      \\
%        0       & \uuu_k^\perp & 0 & 0 & \ldots & r \dot{\uuu}_k^\perp
%      \end{array} \right| 
       &=  \abs{ \left|
      \begin{array}{cc:cccc}
        E_{n,k}   & \uuu_0       & r \dot{\uuu}_0     & 0 & \ldots & 0 \\
        E_{n,k}   & \uuu_1       & 0 & r \dot{\uuu}_1     & \ldots & 0 \\
        \vdots & \vdots & \vdots & \vdots                 & \ddots & \vdots \\
        E_{n,k}   & \uuu_k       & 0 & 0 & \ldots & r \dot{\uuu}_k
      \end{array} \right|} ,
  \end{align}
  where we write the matrix in block notation,
  with $E_{n,k}$ the $n \times k$ matrix with all elements zero
  and ones in the diagonal.
  Similarly, $\uuu_i$ is a column vector of length $n$,
  $r \dot{\uuu}_i$ is an $n \times (n-1)$ matrix,
  and $0$ is the zero matrix of appropriate size,
  which in this case is an $n \times (n-1)$ matrix.
  Like in \cite{ScWe08}, we extract $r$ from $(k+1)(n-1)$ columns,
  and use the fact that transposing the matrix does not affect
  the determinant to get
  \begin{align}
    \left( \frac{J(y,r,\uuu)}{r^{(k+1)(n-1)}} \right)^2  &=
      \left| \begin{array}{cccc}
        E_{k,n}        & E_{k,n}        & \ldots & E_{k,n}        \\
        \uuu_0^T       & \uuu_1^T       & \ldots & \uuu_k^T       \\ \hdashline
        \dot{\uuu}_0^T & 0              & \ldots & 0              \\
        0              & \dot{\uuu}_1^T & \ldots & 0              \\
        \vdots         & \vdots         & \ddots & \vdots         \\
        0              & 0              & \ldots & \dot{\uuu}_k^T 
      \end{array} \right|
    \cdot
      \left| \begin{array}{cc:cccc}
        E_{n,k}   & \uuu_0       & \dot{\uuu}_0       & 0 & \ldots & 0 \\
        E_{n,k}   & \uuu_1       & 0 & \dot{\uuu}_1       & \ldots & 0 \\
        \vdots & \vdots & \vdots & \vdots                 & \ddots & \vdots \\
        E_{n,k}   & \uuu_k       & 0 & 0 & \ldots & \dot{\uuu}_k       \\
      \end{array} \right| .
  \end{align}
  The orthogonality of the matrices $[\uuu_i \dot{\uuu}_i]$ implies that
  $\uuu_i^T \uuu_i = 1$,
  $\dot{\uuu}_i^T \dot{\uuu}_i = E_{n-1,n-1}$,
  whereas $\uuu_i^T \dot{\uuu}_i$ is the zero row vector of length $n-1$, and
  $\dot{\uuu}_i^T \uuu_i$ is the zero column vector of length $n-1$,
  for each $0 \leq i \leq k$.
  We can therefore multiply the matrices and get
  \begin{align}
    \left( \frac{J(y,r,\uuu)}{r^{(k+1)(n-1)}} \right)^2  &=
    \left|
      \begin{array}{cc : ccc}
        (k+1) E_{k,k}   & \sum \uuu_i' & \dot{\uuu}_0' & \ldots & \dot{\uuu}_k'\\
         \sum \uuu_i'^T & k+1          & 0             & \ldots & 0            \\ \hdashline
        \dot{\uuu}_0'^T& 0            & E_{n-1,n-1}   & \ldots & 0            \\
        \vdots         & \vdots       & \vdots        & \ddots & \vdots       \\
        \dot{\uuu}_k'^T& 0            & 0             & \ldots & E_{n-1,n-1}
      \end{array} \right| ,
    \label{eqn:determinant}
  \end{align}
  in which we write $\uuu_i'$ for the vector consisting of the first $k$
  coordinates of $\uuu_i$.
  Similarly, $\dot{\uuu}_i'$ is the $k \times (n-1)$ matrix obtained
  from $\dot{\uuu}_i$ by dropping the bottom $n-k$ rows.
  As written, the $n(k+1) \times n(k+1)$ matrix in \eqref{eqn:determinant}
  is a $(k+3) \times (k+3)$ matrix of blocks, not all of the same size.
  To zero out the last $k+1$ blocks in the first row,
  we subtract the third row times $\dot{\uuu}_0'$,
  the fourth row times $\dot{\uuu}_1'$, and so on.
  The determinant is therefore the product of the determinants of the
  upper left $2 \times 2$ block matrix and the lower right $(k+1) \times (k+1)$
  block matrix, the latter being $1$.
  To further simplify the $2 \times 2$ block matrix,
  we use $[\uuu_i \dot{\uuu}_i] [\uuu_i \dot{\uuu}_i]^T = E_{n,n}$,
  which implies $[\uuu_i' \dot{\uuu}_i'] [\uuu_i' \dot{\uuu}_i']^T = E_{k,k}$,
  and we write the matrix as a product of two matrices:
  \begin{align}
    \left( \frac{J(y,r,\uuu)}{r^{(k+1)(n-1)}} \right)^2
      &= \left| \begin{array}{cc}
           (k+1) E_{k,k} - \sum \dot{\uuu}_i' \dot{\uuu}_i'^T & \sum \uuu_i'  \\
           \sum \uuu_i'^T                                 & k+1 
         \end{array} \right|
    \label{eqn:determinant4} \\
      &= \left| \begin{array}{cc}
           \sum \uuu_i' \uuu_i'^T & \sum \uuu_i'  \\
           \sum \uuu_i'^T         & k+1 
         \end{array} \right|            
       = \left| \left[ \begin{array}{cccc}
                  \uuu_0'  &  \uuu_1'  &  \ldots  &  \uuu_k'  \\  
                  1        &  1        &  \ldots  &  1     
                \end{array} \right]
                \left[ \begin{array}{cc}
                  \uuu_0'^T  &  1       \\  
                  \vdots     &  \vdots  \\
                  \uuu_1'^T  &  1       \\  
                  \uuu_k'^T  &  1 
                \end{array} \right]
         \right| ,
    \label{eqn:determinant5}
  \end{align}
  in which we get from \eqref{eqn:determinant4} to \eqref{eqn:determinant5}
  using $\dot{\uuu}_i' \dot{\uuu}_i'^T = E_{k,k} - \uuu_i' \uuu_i'^T$.
  Finally, the determinant of the vectors $\uuu_i'$ with appended $1$
  is $k!$ times the $k$-dimensional volume of $\uuu'$.
  Hence, $J(y,r,\uuu) = r^{(k+1)(n-1)} k! \Vol{k}{\uuu'}$,
  as claimed in \eqref{eqn:Jacobian}.
  This completes the proof of \eqref{eqn:BP1}.
\eop

\ourparagraph{General case.}
Next we generalize to the case $m \leq k$.
Recall that for a sequence $\xxx$ of $m+1 \leq k+1$ points in $\Rspace^n$,
there is a unique smallest anchored sphere passing through them.
We claim that its center lies inside the orthogonal projection $P$
of the $m$-dimensional affine hull of $\xxx$ onto $\Rspace^k$.
Indeed, orthogonally projecting the center of any anchored sphere
passing through $\xxx$ to $P$ in $\Rspace^k$
we clearly get a point, which is a center of a smaller anchored sphere
still passing through $\xxx$.
The following theorem tells us how to integrate over these smallest
anchored circumscribed spheres.
\begin{theorem}[Anchored Blaschke--Petkantschin Formula]
  \label{thm:BPforGeneralSimplices}
  Let $0 \leq m \leq k \leq n$ and $\alpha = n(m+1)-(k+1)$.
  Then every measurable non-negative function
  $f \colon (\Rspace^n)^{m+1} \to \Rspace$ satisfies
  \begin{align}
    \!\!\!\!\int\displaylimits_{\xxx \in (\Rspace^n)^{m+1}}
          \!\!\!\!\!\!f(\xxx) \diff \xxx
      &=  \int\displaylimits_{y \in \Rspace^k}
          \int\displaylimits_{P \in \LGrass{m}{k}}
          \int\displaylimits_{r \geq 0}
          \int\displaylimits_{\uuu \in (\Sphere)^{m+1}}
            \!\!\!\!\!f(y+r\uuu) r^\alpha [m! \Vol{m}{\uuu'}]^{k-m+1}
            \diff \uuu \diff r \diff P \diff y ,
  \end{align}
  in which $\LGrass{m}{k}$ is the Grassmannian of (linear) $m$-planes
  in $\Rspace^k$, $\uuu'$ is the projection of $\uuu$ to $P$,
  and $\Sphere$ is short for the unit sphere in $P \times \Rspace^{n-k}$.
\end{theorem}
\ourproof
  We use Blaschke--Petkantschin formula twice, first in its standard form.
  For $P \in \LGrass{m}{k}$, we write $P \times \Rspace^{n-k} \in \LGrass{m+n-k}{n}$
  for the $(m+n-k)$-plane whose orthogonal projection to $\Rspace^k$ is $P$.
  The first application of Blaschke--Petkantschin formula integrates over
  all (affine) $m$-planes in $\Rspace^k$,
  spanned by the projections of $\xxx$ to $\Rspace^k$:
  \begin{align}
    \int\displaylimits_{\xxx \in (\Rspace^n)^{m+1}} f(\xxx) \diff \xxx
      &=  \int\displaylimits_{P \in \LGrass{m}{k}}
          \int\displaylimits_{h \in P^\perp}
          \int\displaylimits_{\xxx \in (P \times \Rspace^{n-k})^{m+1}}
            \!\!\!\! f(h+\xxx)
            [m! \Vol{m}{\xxx'}]^{k-m} \diff \xxx \diff h \diff P .
  \end{align}
  For every $m$-plane $P$ in $\Rspace^k$, we consider the vertical
  $(m+n-k)$-plane $P \times \Rspace^{n-k}$ in $\Rspace^n$ and apply
  Lemma \ref{lem:BPforTopdimensionalSimplices} inside it.
  Recalling that $\Sphere$ is the unit sphere in $P \times \Rspace^{n-k}$,
  this gives
  \begin{align}
    \int\displaylimits_{\xxx \in (\Rspace^n)^{m+1}} f(\xxx) \diff \xxx
      &=  \int\displaylimits_{P \in \LGrass{m}{k}}
          \int\displaylimits_{h \in P^\perp}
          \int\displaylimits_{z \in P}
          \int\displaylimits_{r \geq 0}
          \int\displaylimits_{\uuu \in (\Sphere)^{m+1}} f(h+z+r\uuu) 
            r^{(m+n-k-1)(m+1)}                                     \\
      &   ~~~~~~~~~~~~~~~~~~ m! \Vol{m}{\uuu'}
                                 [m! \Vol{m}{r \uuu'}]^{k-m}
            \diff \uuu \diff r \diff z \diff h \diff P .
  \end{align}
  Note that $\Vol{m}{r \uuu'} = r^m \Vol{m}{\uuu'}$, which implies that
  the final power of $r$ is $(m+n-k-1)(m+1) +m(k-m) = \alpha$.
  Finally, we get the claimed relation by setting $y = z + h$ and
  exchanging the integral over $P \in \LGrass{m}{k}$
  with the integral over $y \in \Rspace^k$.
\eop

%%%%%%%%%%%%%%%%%%%%%%%%%%%%%%%%%%%%%%%%%%%%%%%%%%%%%%%%%%%%%%%%%%%%%%%%%%
%%%%%%%%%%%%%%%%%%%%%%%%%%%%%%%%%%%%%%%%%%%%%%%%%%%%%%%%%%%%%%%%%%%%%%%%%%
\section{Expected Number of Intervals}
\label{sec:5}
%%%%%%%%%%%%%%%%%%%%%%%%%%%%%%%%%%%%%%%%%%%%%%%%%%%%%%%%%%%%%%%%%%%%%%%%%%
%%%%%%%%%%%%%%%%%%%%%%%%%%%%%%%%%%%%%%%%%%%%%%%%%%%%%%%%%%%%%%%%%%%%%%%%%%

In this section, we use the anchored Blaschke--Petkantschin formula
of the previous section to compute the expected numbers of intervals
of a weighted Delaunay mosaic in $\Rspace^k$.
We do this for every type and use a weighted Delaunay radius threshold
to get more detailed probabilistic information.
Recall that the weighted mosaic is a random $k$-dimensional slice
of the (unweighted) Poisson--Delaunay mosaic with density
$\density > 0$ in $\Rspace^n$.

\ourparagraph{Slivnyak--Mecke formula.}
To count the type $(\ell, m)$ intervals, we focus our attention
by restricting the center of the weighted Delaunay sphere to a region
$\Region \subseteq \Rspace^k$ and the weighted Delaunay radius to be less than or equal $r_0$.
Any sequence $\xxx = (\xxx_0, \xxx_1, \ldots, \xxx_m)$ of $m+1$
points in $\Rspace^n$ defines such an interval if it
satisfies the following conditions:
\smallskip \begin{enumerate}
  \item[1.] the smallest anchored sphere passing through $\xxx$ is empty,
    writing $\Pempty{\xxx}$ for the probability of this event;
  \item[2.] the center $z$ of this sphere lies in $\Region$,
    writing $\One_\Region (\xxx)$ for the indicator;
  \item[3.] the radius $r$ of this sphere is bounded from above by $r_0$,
    writing $\One_{r_0} (\xxx)$ for the indicator;
  \item[3.] the origin of $\Rspace^k$ sees exactly $m-\ell$ facets of the
    projected $m$-simplex from the outside,
    writing $\One_{m-\ell} (\xxx)$ for the indicator.
\end{enumerate} \smallskip
These are the same conditions as in \cite{ENR16} and \cite{BaEd15}
with the only difference that the sphere is now required to be anchored,
and modulo this remark the proofs are identical.
Combining these conditions with the Slivnyak--Mecke formula,
we get an integral expression for the expected number of type
$(\ell, m)$ intervals, which we partially evaluate
using Theorem \ref{thm:BPforGeneralSimplices}
and Lemma \ref{lem:GammaFunction}:
\begin{align}
  \MoveEqLeft\!\!\!\!\Expected{\ccon{\ell}{m}{k}{n} (r_0)}
    =  \tfrac{1}{(m+1)!}
        \int\displaylimits_{\xxx \in (\Rspace^n)^{m+1}}
          \Pempty{\xxx} \One_\Region (\xxx) \One_{r_0} (\xxx)
          \One_{m-\ell} (\xxx)  \diff \xxx
    \label{eqn:SM1} \\
    &=  \norm{\Region} \norm{\LGrass{m}{k}} \density^{m+1}
          \tfrac{{m!}^{k-m+1}}{(m+1)!}
        \!\!\int\displaylimits_{r \leq r_0}\!\!
          \!e^{- \density \nu_n r^n} r^\alpha \diff r
        \!\!\!\!\!\!\int\displaylimits_{\uuu \in (\Sphere)^{m+1}}\!\!\!\!\!\!\!
          \One_{m-\ell} (\uuu) {\Vol{m}{\uuu'}}^{k-m+1} \diff \uuu
    \label{eqn:SM2} \\
    &=  \norm{\Region} \density^{\frac{k}{n}} \tfrac{{m!}^{k-m}}{m+1}
          \norm{\LGrass{m}{k}}
          \tfrac{\iGama{m+1-\frac{k}{n}}{\density \nu_n r_0^n}}
               {n \nu_n^{m+1-\frac{k}{n}}}
        \!\!\!\int\displaylimits_{\uuu \in (\Sphere)^{m+1}}\!\!\!\!
          \One_{m-\ell} (\uuu) {\Vol{m}{\uuu'}}^{k-m+1} \diff \uuu
    \label{eqn:SM3} \\
    &=  \Ccon{\ell}{m}{k}{n}
          \cdot \tfrac{\iGama{m+1-\frac{k}{n}}{\density \nu_n r_0^n}}
                      { \Gama{m+1-\frac{k}{n}}}
          \cdot \norm{\Region} \density^{\frac{k}{n}} .
    \label{eqn:SM4}
\end{align}
Specifically, we get \eqref{eqn:SM2} by noting
$\Pempty{\xxx} = e^{- \density \nu_n r^n}$,
applying Theorem \ref{thm:BPforGeneralSimplices} to the
right-hand side of \eqref{eqn:SM1},
collapsing the indicators,
using rotational invariance,
and writing $\Sphere$ for the unit sphere in $\Rspace^{m+n-k}$.
We get \eqref{eqn:SM3} from \eqref{eqn:SM2} by applying
Lemma \ref{lem:GammaFunction} with
$j = \alpha + 1 = n(m+1)-k$, $c = \density \nu_n$, $p = n$, $t_0 = r_0$,
which asserts that the integral over the radius evaluates to
the fraction involving the incomplete Gamma function.
%%$\iGama{k/p}{c t_0^p} / (p c^{k/p})$.
  %% = \iGama{m+1-k/n}{n \density^{m+1-k/n} \nu_n^{m+1-k/n}}$.
Finally, we get \eqref{eqn:SM4} by defining the constant
\begin{align}
  \Ccon{\ell}{m}{k}{n}  &= 
    \tfrac{{m!}^{k-m} \norm{\LGrass{m}{k}} \Gama{m+1-\frac{k}{n}}} 
          {(m+1) n \nu_n^{m+1-\frac{k}{n}}}
    \int\displaylimits_{\uuu \in (\Sphere)^{m+1}}
          \One_{m-\ell} (\uuu) {\Vol{m}{\uuu'}}^{k-m+1} \diff \uuu .
    \label{eqn:constantlm}
\end{align}
As a sanity check, we set $\ell = m = 0$ and $k = 1$, and get
$\Ccon{0}{0}{1}{n} = \sigma_{n-1} \Gama{1-1/n} / (n \nu_n^{1-1/n})$
because $\Sphere \subseteq \Rspace^{n-1}$ has volume $\sigma_{n-1}$,
and we have $\One_{0} (\uuu_0) = 1$ and $\Vol{0}{\uuu_0} = 1$
for all points $\uuu_0 \in \Sphere$.
This agrees with \eqref{eqn:constant00} in Section \ref{sec:2}.

\ourparagraph{Simplices in the weighted Delaunay mosaic.}
Since this constant in \eqref{eqn:constantlm} does not depend on $r_0$,
we deduce that the weighted Delaunay radius of a typical type $(\ell, m)$ interval is
Gamma distributed.
The weighted Delaunay radius of a typical $j$-simplex in the weighted Poisson--Delaunay mosaic
therefore follows a linear combination of Gamma distributions.
Indeed, we get the total number of $j$-simplices as
$\dcon{j}{k}{n} = \sum_{\ell=0}^j \sum_{m=j}^k
                  \binom{m-\ell}{m-j} \ccon{\ell}{m}{k}{n}$; see \cite{ENR16}.
The same relation holds if we limit the simplices to weighted Delaunay radius at most $r_0$,
and also if we replace the simplex counts by the constants
$\Ccon{\ell}{m}{k}{n}$ and the analogously defined $\Dcon{j}{k}{n}$.
Before continuing, we consider the top-dimensional case, $j = k$,
in which $\Dcon{k}{k}{n} = \sum_{\ell=0}^k \Ccon{\ell}{k}{k}{n}$.
Taking the sum eliminates the indicator function in \eqref{eqn:constantlm},
and we get
\begin{align}
  \Dcon{k}{k}{n} 
    &=  \frac{\Gama{k+1-\frac{k}{n}}}{(k+1) n \nu_n^{k+1-\frac{k}{n}}}
        \int\displaylimits_{\uuu \in (\Sspace^{n-1})^{k+1}}
          \Vol{k}{\uuu'} \diff \uuu .
  \label{eqn:ksimplices1}
\end{align}
We can compare this with the expression for the number of Voronoi
vertices by M{\o}ller \cite{Mol89} using Crofton formula \cite[Chapter 6]{Had57};
see also \cite[Theorem 10.2.4]{ScWe08}.
By duality, the number of vertices in the weighted Voronoi tessellation
is the number of top-dimensional simplices in the
weighted Delaunay mosaic.
Each vertex is the intersection of an $(n-k)$-dimensional Voronoi polyhedron
with the $k$-plane, and if we know the expected number of intersections,
then we also know the integral, over all $k$-planes.
Crofton formula applies and gives the $(n-k)$-dimensional volume
of the $(n-k)$-skeleton of the Voronoi tessellation
as $\sigma_n / (2 \norm{\LGrass{k}{n}} \nu_{n-1})$ times the
mentioned integral.
It turns out that the expected volume is not so difficult
to compute otherwise \cite{Mol89},
so we can turn the argument around and deduce the expected number
of vertices from the expected $(n-k)$-dimensional volume.
This gives
\begin{align}
  \Dcon{k}{k}{n} 
    &=  \frac{\sigma_1 \sigma_{n+1}}{\sigma_{k+1} \sigma_{n-k+1}}
        \frac{2^{k+1} \pi^{k/2}}{n (k+1)!}
        \frac{\Gama{\frac{kn+n-k+1}{2}}}{\Gama{\frac{kn+n-k}{2}}}
        \frac{\Gama{\frac{n+2}{2}}^{k+1-\frac{k}{n}}}{\Gama{\frac{n+1}{2}}^k}
        \frac{\Gama{k+1-\frac{k}{n}}}{\Gama{\frac{n-k+1}{2}}} .
  \label{eqn:ksimplices2}
\end{align}
Comparing \eqref{eqn:ksimplices2} with \eqref{eqn:ksimplices1},
we get an explicit expression for the expected $k$-dimensional
volume of the projection of a random $k$-simplex inscribed in $\Sspace^{n-1}$.

\section{Computations}
\label{sec:6}
%%%%%%%%%%%%%%%%%%%%%%%%%%%%%%%%%%%%%%%%%%%%%%%%%%%%%%%%%%%%%%%%%%%%%%%%%%
%%%%%%%%%%%%%%%%%%%%%%%%%%%%%%%%%%%%%%%%%%%%%%%%%%%%%%%%%%%%%%%%%%%%%%%%%%

We now return to \eqref{eqn:constantlm} and note that the integral
on the right-hand side is $\sigma_{m+n-k}^{m+1}$ times the expected
value of the random variable
\begin{align}
  \Random{\ell}{m}{k}{n}  &=  \One_{m-\ell} (\uuu) \Vol{m}{\uuu'}^{k-m+1} ,
  \label{eqn:random}
\end{align}
where $\uuu$ is a sequence of $m+1$ random points uniformly and
independently distributed on the unit sphere in $\Rspace^{m+n-k}$,
and $\uuu'$ is the corresponding sequence of points projected to
$\Rspace^m \hookrightarrow \Rspace^{m+n-k}$.
Our goal is to compute $\Expected{\Random{\ell}{m}{k}{n}}$
in some special cases.
Instead of working with the original points, we prefer to study
their projections to $\Rspace^m$,
but the distribution of the $m+1$ points in $\Rspace^m$ has yet to
be determined.
If the upper bound is a vertex or an edge, then we find
explicit expressions of the expected number of intervals.

\ourparagraph{Critical vertices.}
For $m = 0$, we count intervals of type $(0,0)$ or, equivalently,
critical vertices.
Since $\Random{0}{0}{k}{n} = 1$, for all $k \leq n$, we get
\begin{align}
  \Ccon{0}{0}{k}{n}  &=  \sigma_{n-k}
                         \tfrac{ \Gama{1-\frac{k}{n}} } { n \nu_n^{1-k/n} }
  \label{eqn:criticalvertices}
\end{align}
from \eqref{eqn:constantlm}.
Accordingly, the expected number of critical vertices in $\Region$
with weighted Delaunay radius at most $r_0$ is $\Ccon{0}{0}{k}{n}$
times the normalized incomplete Gamma function
times $\norm{\Region} \density^{k/n}$;
compare with \eqref{eqn:critical-vertices-1}
and \eqref{eqn:critical-vertices-2} in Section \ref{sec:2}.

\ourparagraph{Vertex-edge pairs.}
Next we count the intervals of type $(0,1)$ or, equivalently,
the regular vertex-edge pairs.
For this, we need the expectation of $\Random{0}{1}{k}{n}$:
picking two random points on the unit sphere in $\Rspace^{n-k+1}$
and projecting them to $\Rspace^1 \hookrightarrow \Rspace^{n-k+1}$,
this is the expectation when we get the $k$-th power of the distance
between the projected points, if they lie on the same side of the origin,
and we get $0$, otherwise.
Writing $\uuu_0', \uuu_1' \in [-1,1]$ for the projected points
and $x = |\uuu_0'|$, $y = |\uuu_1'|$ for their absolute values,
we note that the signs and magnitudes are independent.
It follows that we get zero with probability $\tfrac{1}{2}$,
so the desired expectation is
\begin{align}
  \Expected{\Random{0}{1}{k}{n}}
    &=  \tfrac{1}{2} \Expected{ |x-y|^k }
     =  \Expected{ (x-y)^k \One_{x > y} } .
\end{align}
We can therefore restrict our attention to the half of the unit sphere
that projects to $[0,1]$.
To integrate over this hemisphere, we use that $x^2$ and $y^2$ are independent
Beta-distributed random variables; see Appendix \ref{app:A}.
Setting $a = x^2$ and $b = y^2$, we have
\begin{align}
  \Expected{\Random{0}{1}{k}{n}}
    &=  \frac{1}{\Beta{\tfrac{n-k}{2}}{\tfrac{1}{2}}^2}
        \int\displaylimits_{a=0}^1
        \int\displaylimits_{b=0}^a [\sqrt{a}-\sqrt{b}]^k
          a^{-\frac{1}{2}} (1-a)^{\frac{n-k-2}{2}}
          b^{-\frac{1}{2}} (1-b)^{\frac{n-k-2}{2}} \diff a \diff b 
    \label{eqn:vertexedgepairs-1} \\
    &=  \frac{4}{\Beta{\tfrac{n-k}{2}}{\tfrac{1}{2}}^2}
        \int\displaylimits_{x=0}^1
        \int\displaylimits_{y=0}^x [x-y]^k
          (1-x^2)^{\frac{n-k-2}{2}}
          (1-y^2)^{\frac{n-k-2}{2}} \diff x \diff y 
    \label{eqn:vertexedgepairs-2} \\
    &=  \frac{ \Gama{k+1} \Gama{\tfrac{n-k+1}{2}}^2 }
             { 2^k \sqrt{\pi} \Gama{\tfrac{n-k}{2}} } \cdot
        \HyperReg{3}{2} \left( \tfrac{1}{2}, 1, \tfrac{k-n+2}{2};
                      \, \tfrac{k+3}{2}, \tfrac{n+2}{2}; \, 1 \right) ,
    \label{eqn:vertexedgepairs-3}
\end{align}
in which $\HyperReg{3}{2}$ is the regularized hypergeometric function
considered in Appendix \ref{app:A} and we use the Mathematica software
to get from \eqref{eqn:vertexedgepairs-2} to \eqref{eqn:vertexedgepairs-3}.
As mentioned at the end of this appendix,
$\tfrac{k+3}{2} + \tfrac{n+2}{2} > \tfrac{1}{2} + 1 + \tfrac{k-n+2}{2}$
is a sufficient condition for the convergence of the infinite sum
that defines the value of the regularized hypergeometric function.
This is equivalent to $n > 0$, which is always satisfied.
Plugging \eqref{eqn:vertexedgepairs-3} into \eqref{eqn:constantlm},
we get an expression for the corresponding constant:
\begin{align}
  \Ccon{0}{1}{k}{n}
    &=  \frac{\sigma_{n-k+1}^2 \sigma_k \Gama{2-\frac{k}{n}}}
             {4 n \nu_n^{2-k/n}}
        \frac{ \Gama{k+1} \Gama{\tfrac{n-k+1}{2}}^2 }
             { 2^k \sqrt{\pi} \Gama{\tfrac{n-k}{2}} } \cdot
        \HyperReg{3}{2} \left( \tfrac{1}{2}, 1, \tfrac{k-n+2}{2};
                      \, \tfrac{k+3}{2}, \tfrac{n+2}{2}; \, 1 \right) .
    \label{eqn:vertexedgepairs-4}
\end{align}

\ourparagraph{Critical edges.}
Next we count the intervals of type $(1,1)$ or, equivalently,
the critical edges.
Here the expectation of $\Random{1}{1}{k}{n}$ is relevant:
picking two points on the unit sphere in $\Rspace^{n-k+1}$ and projecting
them to $\Rspace^1 \hookrightarrow \Rspace^{n-k+1}$,
this is the expectation in which we get the $k$-th power of the distance
between the projected points, if they lie on opposite sides of the origin,
and we get $0$, otherwise.
Using again that the signs and magnitude of the projected points are
independent, we note that this expectation is
$\Expected{\Random{1}{1}{k}{n}}  =  \tfrac{1}{2} \Expected{(x+y)^k}$.
Setting $a = x^2$, $b = y^2$, and integrating as before, we get
\begin{align}
  \MoveEqLeft\Expected{\Random{1}{1}{k}{n}}
    =  \frac{1}{\Beta{\tfrac{n-k}{2}}{\tfrac{1}{2}}^2}
        \int\displaylimits_{a=0}^1
        \int\displaylimits_{b=0}^1
          \left[ \sqrt{a} + \sqrt{b} \right]^k
                                  a^{-\frac{1}{2}} (1-a)^{\frac{n-k-2}{2}}
                                  b^{-\frac{1}{2}} (1-b)^{\frac{n-k-2}{2}}
          \diff a \diff b
  \label{eqn:criticaledges-1} \\
    &=  \frac{1}{\Beta{\tfrac{n-k}{2}}{\tfrac{1}{2}}^2}
        \int\displaylimits_{a=0}^1
        \int\displaylimits_{b=0}^1
          \sum\displaylimits_{i=0}^k \binom{k}{i}
                                  a^{\frac{i-1}{2}} b^{\frac{k-i-1}{2}}
                                  (1-a)^{\frac{n-k-2}{2}} (1-b)^{\frac{n-k-2}{2}}
          \diff a \diff b
  \label{eqn:criticaledges-2} \\
    &=  \frac{1}{\Beta{\tfrac{n-k}{2}}{\tfrac{1}{2}}^2}
          \sum\displaylimits_{i=0}^k \binom{k}{i}
                                  \Beta{\tfrac{n-k}{2}}{\tfrac{i+1}{2}}
                                  \Beta{\tfrac{n-k}{2}}{\tfrac{k-i+1}{2}} .
  \label{eqn:criticaledges-3}
\end{align}
Plugging \eqref{eqn:criticaledges-3} into \eqref{eqn:constantlm}, we
get the expression for the corresponding constant:
\begin{align}
  \MoveEqLeft\Ccon{1}{1}{k}{n}  =  \frac{\sigma_{n-k+1}^2 \sigma_k \Gama{2-\frac{k}{n}}}
                              {8 n \nu_n^{2-k/n} \Beta{\frac{n-k}{2}}{\frac{1}{2}}^2}
                         \sum\displaylimits_{i=0}^k \binom{k}{i}
                                  \Beta{\tfrac{n-k}{2}}{\tfrac{i+1}{2}}
                                  \Beta{\tfrac{n-k}{2}}{\tfrac{k-i+1}{2}} .
  \label{eqn:criticaledges-4}
\end{align}

\ourparagraph{Constants in low dimensions.}
The authors have checked the $k$-dimensional formulas 
against the $1$-dimensional formulas in Section \ref{sec:2},
both symbolically and numerically.
In $k=2$ dimensions, the formulas provide sufficient information
to compute all constants governing the expectations of the six types
of intervals.
We get three constants from \eqref{eqn:criticalvertices},
\eqref{eqn:vertexedgepairs-4}, \eqref{eqn:criticaledges-4}:
\begin{align}
  \Ccon{0}{0}{2}{n}  &=  \frac{\sigma_{n-2} \Gama{1-\frac{2}{n}}}
                              {n \nu_n^{1-2/n}} , 
  \label{eqn:Constantk2-00} \\
  \Ccon{0}{1}{2}{n}  &=  \frac{ \sigma_{n-1}^2 \sqrt{\pi} \Gama{2-\frac{2}{n}}}
                              { 4 n \nu_n^{2-2/n}} 
                         \frac{ \Gama{\tfrac{n-1}{2}}^2 }
                              { \Gama{\tfrac{n-2}{2}} } \cdot
                         \HyperReg{3}{2} \left( \tfrac{1}{2}, 1, \tfrac{4-n}{2};
                                       \, \tfrac{5}{2}, \tfrac{n+2}{2}; \, 1 \right) ,
  \label{eqn:Constantk2-01} \\
  \Ccon{1}{1}{2}{n}  &=  \frac{\sigma_{n-1}^2 \Gama{2-\frac{2}{n}} \pi}
                              {2 n \nu_n^{2-2/n}} \cdot
                         \left[ \frac{1}{n-1} +
                                \frac{\Gama{\frac{n-1}{2}}^2}
                                     {\pi \Gama{\frac{n}{2}}^2} \right] .
  \label{eqn:Constantk2-11}
\end{align}
The critical simplices satisfy the Euler relation \cite{For98}:
$\Ccon{0}{0}{2}{n} - \Ccon{1}{1}{2}{n} + \Ccon{2}{2}{2}{n} = 0$,
which gives us the constant for the critical triangles.
We get another linear relation from the fact that in the plane the number
of triangles is twice the number of vertices \cite[page 458, Theorem 10.1.2]{ScWe08}:
$\Ccon{0}{2}{2}{n} + \Ccon{1}{2}{2}{n} + \Ccon{2}{2}{2}{n} =
  2( \Ccon{0}{0}{2}{n} + \Ccon{0}{1}{2}{n} + \Ccon{0}{2}{2}{n})$.
Finally, we get a relation for the number of weighted Delaunay triangles
from \eqref{eqn:ksimplices2}, which we restate for $k=2$:
\begin{align}
  \Dcon{2}{2}{n} 
    &=  \frac{2 \sigma_{n+1}}{3 n \sigma_{n-1}}
        \frac{\Gama{\frac{3n-1}{2}}}{\Gama{\frac{3n-2}{2}}}
        \frac{\Gama{\frac{n+2}{2}}^{3-\frac{2}{n}}}{\Gama{\frac{n+1}{2}}^2}
        \frac{\Gama{3-\frac{2}{n}}}{\Gama{\frac{n-1}{2}}} .
\end{align}
Combining $\Ccon{0}{2}{2}{n} + \Ccon{1}{2}{2}{n} + \Ccon{2}{2}{2}{n} = \Dcon{2}{2}{n}$
with the two linear relations mentioned above, we get
\begin{align}
  \Ccon{0}{2}{2}{n}  &=  - \Ccon{0}{0}{2}{n} - \Ccon{0}{1}{2}{n}
                         + \tfrac{1}{2} \Dcon{2}{2}{n} ,                     \\
  \Ccon{1}{2}{2}{n}  &=    \Ccon{0}{0}{2}{n} + \Ccon{0}{1}{2}{n}
                         - \Ccon{2}{2}{2}{n} + \tfrac{1}{2} \Dcon{2}{2}{n} , \\
  \Ccon{2}{2}{2}{n}  &=  - \Ccon{0}{0}{2}{n} + \Ccon{1}{1}{2}{n} .
\end{align}
For small values of $n$, the constants are approximated in Table \ref{tbl:Constants2D}.
\begin{table}[hbt]
  \centering
  \small{ \begin{tabular}{r||rrrr rrrr rrrr}
    & \multicolumn{1}{c}{$n=3$}& \multicolumn{1}{c}{$4$} & \multicolumn{1}{c}{$5$}   
    & \multicolumn{1}{c}{$6$}  & \multicolumn{1}{c}{$7$} & \multicolumn{1}{c}{$8$}
    & \multicolumn{1}{c}{$9$}  & \multicolumn{1}{c}{$10$}& \multicolumn{1}{c}{$\ldots$}
    & \multicolumn{1}{c}{$20$} & \multicolumn{1}{c}{$\ldots$}& \multicolumn{1}{c}{$1000$}
                                                                  \\ \hline \hline
    $\Ccon{0}{0}{2}{n}$
      & $1.11$ & $1.25$ & $1.38$ & $1.49$ & $1.58$ & $1.66$ & $1.73$ & $1.79$
                                 &$\ldots$& $2.12$ &$\ldots$& $2.69$  \\
    $\Ccon{0}{1}{2}{n}$
      & $0.26$ & $0.42$ & $0.54$ & $0.63$ & $0.71$ & $0.77$ & $0.82$ & $0.86$
                                 &$\ldots$& $1.12$ &$\ldots$& $1.54$  \\
    $\Ccon{0}{2}{2}{n}$
      & $0.09$ & $0.15$ & $0.21$ & $0.25$ & $0.28$ & $0.31$ & $0.33$ & $0.35$
                                 &$\ldots$& $0.47$ &$\ldots$& $0.65$  \\
    $\Ccon{1}{1}{2}{n}$
      & $2.47$ & $2.92$ & $3.30$ & $3.61$ & $3.87$ & $4.09$ & $4.28$ & $4.44$
                                 &$\ldots$& $5.37$ &$\ldots$& $6.92$  \\
    $\Ccon{1}{2}{2}{n}$
      & $1.46$ & $1.83$ & $2.13$ & $2.37$ & $2.57$ & $2.74$ & $2.89$ & $3.01$
                                 &$\ldots$& $3.72$ &$\ldots$& $4.88$  \\
    $\Ccon{2}{2}{2}{n}$
      & $1.37$ & $1.67$ & $1.92$ & $2.12$ & $2.29$ & $2.43$ & $2.55$ & $2.66$
                                 &$\ldots$& $3.25$ &$\ldots$& $4.23$  \\ \hline
    $\Dcon{0}{2}{n}$
      & $1.46$ & $1.83$ & $2.13$ & $2.37$ & $2.57$ & $2.74$ & $2.89$ & $3.01$
                                 &$\ldots$& $3.72$ &$\ldots$& $4.88$  \\
    $\Dcon{1}{2}{n}$
      & $4.37$ & $5.48$ & $6.38$ & $7.10$ & $7.71$ & $8.22$ & $8.66$ & $9.03$
                                 &$\ldots$&$11.16$ &$\ldots$&$14.65$  \\
    $\Dcon{2}{2}{n}$
      & $2.92$ & $3.66$ & $4.25$ & $4.74$ & $5.14$ & $5.48$ & $5.77$ & $6.02$
                                 &$\ldots$& $7.44$ &$\ldots$& $9.77$ 
  \end{tabular} }
  \caption{The rounded constants in the expressions of the expected number of
    intervals and simplices of a $2$-dimensional weighted Delaunay mosaic obtained
    from a Poisson point process in $n$ dimensions.}
  \label{tbl:Constants2D}
\end{table}

%%%%%%%%%%%%%%%%%%%%%%%%%%%%%%%%%%%%%%%%%%%%%%%%%%%%%%%%%%%%%%%%%%%%%%%%%%
%%%%%%%%%%%%%%%%%%%%%%%%%%%%%%%%%%%%%%%%%%%%%%%%%%%%%%%%%%%%%%%%%%%%%%%%%%
\section{Discussion}
\label{sec:7}
%%%%%%%%%%%%%%%%%%%%%%%%%%%%%%%%%%%%%%%%%%%%%%%%%%%%%%%%%%%%%%%%%%%%%%%%%%
%%%%%%%%%%%%%%%%%%%%%%%%%%%%%%%%%%%%%%%%%%%%%%%%%%%%%%%%%%%%%%%%%%%%%%%%%%

The main result of this paper is the stochastic analysis of the radius
function of a weighted Poisson--Delaunay mosaic.
As a consequence, we get formulas for the expected number of simplices
in weighted Poisson-Delaunay mosaics; compare with \cite{Lau07,LaZu08}.
The main technical steps leading up to this result are a
new Blaschke--Petkantschin formula for spheres,
stated as Theorem \ref{thm:BPforGeneralSimplices},
and the discrete Morse theory approach recently introduced in \cite{ENR16}.
There are a number of open questions that remain:
\begin{itemize}\denselist
  \item  We have explicit expressions for the constants in the expected
    number of intervals of all types for dimension $k \leq 2$.
    To go beyond two dimensions, Wendel's method of reflecting vertices
    of a simplex through the origin \cite{Wen62}  should be useful.
    Short of getting precise formulas, can we say something
    about the asymptotic behavior of the
    constants, as $k$ and $n$ go to infinity?
  \item  The connection to Crofton formula and the volumes of Voronoi
    skeleta has been mentioned in Section \ref{sec:5}.
    Are there further connections that relate such volumes with simplices
    of dimension strictly less than $k$, or with subsets of simplices
    limited to radii at most $r_0$?
  \item  The slice construction implies a repulsive force among the vertices:
    the vertices of the weighted Poisson--Delaunay mosaic
    are more evenly spread than a Poisson point process.
    For fixed $k$, the repulsion gets stronger with increasing $n$.
    It would be interesting to study this repulsive force and
    its consequences analytically.
\end{itemize}

%\newpage
%%%%%%%%%%%%%%%%%%%%%%%%%%%
\bibliographystyle{acm}
\bibliography{weighted}

\appendix
%%%%%%%%%%%%%%%%%%%%%%%%%%%%%%%%%%%%%%%%%%%%%%%%%%%%%%%%%%%%%%%%%%%%%%%%%%
%%%%%%%%%%%%%%%%%%%%%%%%%%%%%%%%%%%%%%%%%%%%%%%%%%%%%%%%%%%%%%%%%%%%%%%%%%
\section{On Special Functions}
\label{app:A}
%%%%%%%%%%%%%%%%%%%%%%%%%%%%%%%%%%%%%%%%%%%%%%%%%%%%%%%%%%%%%%%%%%%%%%%%%%
%%%%%%%%%%%%%%%%%%%%%%%%%%%%%%%%%%%%%%%%%%%%%%%%%%%%%%%%%%%%%%%%%%%%%%%%%%

In this appendix, we define and discuss three types of special functions
used in the main body of this paper:
Gamma functions, Beta functions, and hypergeometric functions.

\ourparagraph{Gamma functions.}
We recall that the \emph{lower-incomplete Gamma function} takes two parameters,
$j$ and $t_0 \geq 0$, and is defined by
\begin{align}
  \iGama{j}{t_0}  &=  \int_{t=0}^{t_0} t^{j-1} e^{-t} \diff t .
  \label{eqn:iGama}
\end{align}
The corresponding complete \emph{Gamma function} is
$\Gama{j} = \iGama{j}{\infty}$.
An important relation for Gamma functions is $\Gama{j+1} = j \Gama{j}$,
which holds for any real $j$ that is not a non-positive integer.
We often use the ratio, $\iGama{j}{t_0} / \Gama{j}$,
which is the density of a probability distribution
and called the \emph{Gamma distribution} with parameter $j$.
We prove a technical lemma about incomplete Gamma functions,
which is repeatedly used in the main body of this paper.
\begin{lemma}[Gamma Function]
  \label{lem:GammaFunction}
  Let $c, p, j, t_0 \in \Rspace$ with $p \neq 0$ and $t_0 > 0$.
  Then
  \begin{align}
    \int_{t=0}^{t_0} t^{j-1} e^{-c t^p} \diff t
      &=  \frac{ \iGama{\frac{j}{p}}{c t_0^p} }
               { p c^{j/p} } .
    \label{eqn:A}
  \end{align}
\end{lemma}
\ourproof
  We rewrite the numerator of the right-hand side of \eqref{eqn:A}
  using the definition of the right-incomplete Gamma function
  \eqref{eqn:iGama} and substituting $u = c t^p$ and
  $\diff u = c p t^{p-1} \diff t$:
  \begin{align}
    \iGama{\frac{j}{p}}{c t_0^p} 
      &=  \int_{u=0}^{c t_0^p}
            u^{\frac{j}{p}-1} e^{-u} \diff u             \\
      &=  \int_{t=0}^{t_0}
            (c t^p)^{\frac{j}{p}-1} e^{-c t^p} cp t^{p-1} \diff t \\
      &=  \int_{t=0}^{t_0} p c^{\frac{j}{p}}
            t^{j-1} e^{-c t^p} \diff t .
  \end{align}
  Dividing by $p c^{j/p}$ gives the claimed relation.
\eop

\ourparagraph{Beta functions.}
Given real numbers $a$, $b$, and $0 \leq t_0 \leq 1$,
the \emph{incomplete Beta function} is defined by
\begin{align}
  \iBeta{t_0}{a}{b}  &=  \int_{t=0}^{t_0} t^{a-1} (1-t)^{b-1} \diff t ,
  \label{eqn:Betadefinition}
\end{align}
and the complete \emph{Beta function} is $\Beta{a}{b} = \iBeta{1}{a}{b}$,
which can be expressed in terms of complete
Gamma functions:  $\Beta{a}{b} = \Gama{a} \Gama{b} / \Gama{a+b}$.

The Beta functions can be used to integrate over the projection of
a sphere in $\Rspace^n$ to a linear subspace $\Rspace^k \hookrightarrow \Rspace^n$,
as we now explain.
Assuming $\Rspace^k$ is spanned by the first $k$ coordinate vectors of $\Rspace^n$,
the projection of a point means dropping coordinates $k+1$ to $n$.
Suppose now that we pick a point $x = (x_1, x_2, \ldots, x_n)$
uniformly on $\Sspace^{n-1}$ by normalizing a vector of $n$ normally
distributed random variables:
$X_i \sim \Normal{0}{1}$ for $1 \leq i \leq n$
and $x_j = X_j / \left( \sum_{i=1}^n X_i^2 \right)^{1/2}$ for $1 \leq j \leq n$.
Its projection to $\Rspace^k$ is $x' = (x_1, \ldots, x_k, 0, \ldots, 0)$,
and the squared distance from the origin is
$\norm{x'}^2 = \left( \sum_{i=1}^k x_i^2 \right)
             / \left( \sum_{i=1}^n x_i^2 \right)$.
It can be written as $r^2 = X / (X+Y)$,
in which $X$ and $Y$ are $\chi^2$-distributed independent random variables
with $k$ and $n-k$ degrees of freedom, respectively.
This implies that $r^2 \sim \Beta{\tfrac{k}{n}}{\tfrac{n-k}{n}}$
\cite[Section 4.2]{Wal96}.
Consider for example the case $k = n-1$.
Integrating in $\Rspace^k$ over all points with distance at most $r_0$ from the origin
is the same as integrating over two spherical caps of $\Sspace^{n-1}$,
namely the cap around the north-pole bounded by $(n-2)$-spheres of radius $r_0$,
and a similar cap around the south-pole.
To compute the volume of a single such cap, we set $t_0 = r_0^2$
and integrate the incomplete Beta function:
\begin{align}
  \Vol{n-1}{r_0} &=  \frac{\sigma_n}{2 \Beta{\tfrac{n-1}{2}}{\tfrac{1}{2}}}
                     \int_{t=0}^{t_0} t^{\frac{n-1}{2}-1}
                                  (1-t)^{\frac{1}{2}-1} \diff t 
                 =   \frac{\iBeta{t_0}{\tfrac{n-1}{2}}{\tfrac{1}{2}}}
                          {2   \Beta{\tfrac{n-1}{2}}{\tfrac{1}{2}}}.
\end{align}
Similarly, we can integrate over a ball in a $k$-dimensional projection
and get the volume of the preimage, which is a solid torus
inside the $(n-1)$-sphere.

\ourparagraph{Hypergeometric functions.}
The family of \emph{hypergeometric functions} takes
$p+q$ parameters and one argument and can be defined as a sum of products
of Gamma functions,
while the \emph{regularized} version of this function is obtained by normalizing
by the product of $\Gama{b_i}$:
\begin{align}
  \Hyper{p}{q} \left( a_1, \ldots, a_p; b_1, \ldots, b_q; z \right)
    &=  \sum_{j=0}^\infty \left[ \prod_{i=1}^p \frac{\Gama{j+a_i}}{\Gama{a_i}} \right]
                          \left[ \prod_{i=1}^q \frac{\Gama{b_i}}{\Gama{j+b_i}} \right]
                          \frac{z^j}{j!} , \\
  \HyperReg{p}{q} \left( a_1, \ldots, a_p; b_1, \ldots, b_q; z \right)
    &=  \Hyper{p}{q} \left( a_1, \ldots, a_p; b_1, \ldots, b_q; z \right)
          / \prod_{i=1}^q \Gama{b_i}      \\
    &=  \sum_{j=0}^\infty \left[ \prod_{i=1}^p \frac{\Gama{j+a_i}}{\Gama{a_i}} \right]
                          \left[ \prod_{i=1}^q \frac{1}{\Gama{j+b_i}} \right]
                          \frac{z^j}{j!} .
\end{align}
We are interested in the type $p=3$ and $q=2$.
Here convergence of the infinite sum depends on the values of the parameters.
We always have convergence for $|z| < 1$, and if $z = 1$, a sufficient condition for convergence is $b_1 + b_2 > a_1 + a_2 + a_3$
\cite{OLBC10}.

\end{document}